\documentclass[preprint,12pt]{elsarticle}

\makeatletter
\def\ps@pprintTitle{%
 \let\@oddhead\@empty
 \let\@evenhead\@empty
 \def\@oddfoot{}%
 \let\@evenfoot\@oddfoot}
\makeatother

\usepackage[top=2.5cm,bottom=2.5cm,left=3cm,right=3cm]{geometry}

\usepackage{amsmath,amssymb}
\usepackage{mathrsfs}

\usepackage{listings}

\usepackage{graphicx}
\usepackage{subfig}
\usepackage{algorithm}
\usepackage[noend]{algpseudocode}
\usepackage{enumitem}
\usepackage[hidelinks]{hyperref}
\usepackage{makecell}

\usepackage{pgfplots}
\usepackage{array}

\usepackage{colortbl,array}
\usepackage{multirow,bigdelim}
\usepackage{booktabs}

\makeatletter
\newcommand{\shorteq}{%
  \settowidth{\@tempdima}{-}
  \resizebox{\@tempdima}{\height}{=}%
}
\makeatother

\newcolumntype{H}{>{\setbox0=\hbox\bgroup}c<{\egroup}@{}}

\usepackage{alltt}
\newcommand{\bX}{{\bf X}}

\newcommand{\nullmat}{{\bf S}}

\newcommand{\rel}{\tau}
\newcommand{\brel}{\mathbf{c}}

\newcommand{\mX}{{\bf x}}
\newcommand{\tX}{{\mathcal X}}
\newcommand{\X}{\mathbf{X}}
\newcommand{\mY}{{\bf y}}

\newcommand{\mZ}{{\bf z}}
\newcommand{\tZ}{{\mathcal Z}}
\newcommand{\vv}{{\bf v}}

\usepackage{caption}

\newcommand{\multipleX}{\mathbf{X}}

\newcommand{\mH}{{\bf H}}
\newcommand{\tH}{\mathcal{H}}

\newcommand{\I}{\mathbf{I}}
\newcommand{\bS}{\mathbf{S}}
\newcommand{\M}{\mathcal{M}}
\newcommand{\OP}{\mathsf{P}}

\newcommand{\bsize}{b}
\newcommand{\R}{\mathbf{R}}
\newcommand{\PP}{\mathbf{P}}
\newcommand{\C}{\mathbf{C}}
\newcommand{\Lam}{\mathbf{\Lambda}}
\newcommand{\Pvec}{\mathbf{p}}
\newcommand{\Rvec}{\mathbf{r}}
\newcommand{\Prec}{\mathbf{B}^{-1}}
\newcommand{\Retr}{\mathcal{T}}
\newcommand{\tens}[1]{\mathcal{#1}}
\newcommand{\tensel}[1]{\mathcal{#1}}

\newcommand{\energy}{\epsilon}

\newcommand{\basis}{\mathbf{V}}

\newcommand{\Rnd}{\mathbb{R}^{n^d}}

\newcommand{\RQ}{\mathfrak{R}}

\newcommand{\diag}{\text{diag}}
\DeclareMathOperator{\Tr}{Tr}
\newcommand{\rank}{{r}}

\newcommand{\Mult}{\mathbf{\Lambda}}
\newcommand{\mult}{\lambda}
\newcommand{\Aoo}{\X^\intercal \mH \X}
\newcommand{\Aol}{\X^\intercal \mH \basis}
\newcommand{\Alo}{\basis^\intercal \mH \X}
\newcommand{\All}{\basis^\intercal \mH \basis}
\newcommand{\Goo}{\X^\intercal \X}
\newcommand{\Gol}{\X^\intercal \basis}
\newcommand{\Glo}{\basis^\intercal \X}
\newcommand{\Gll}{\basis^\intercal \basis}

\usepackage[normalem]{ulem}

\usepackage{amsthm}

\usepackage{mathtools}
\DeclarePairedDelimiter\abs{\lvert}{\rvert}
\DeclarePairedDelimiter\norm{\lVert}{\rVert}

\makeatletter
\let\oldabs\abs
\def\abs{\@ifstar{\oldabs}{\oldabs*}}

\let\oldnorm\norm
\def\norm{\@ifstar{\oldnorm}{\oldnorm*}}
\makeatother

\DeclareMathOperator*{\argmin}{arg\,min}
\DeclareMathOperator*{\argmax}{arg\,max}

\begin{document}
\begin{frontmatter}
\title{Low-rank Riemannian eigensolver for high-dimensional Hamiltonians}

\author[add1,fn1]{Maxim Rakhuba\corref{cor1}}
\ead{maksim.rakhuba@sam.math.ethz.ch}
\fntext[fn1]{Work by MR on this project was performed while he was a junior research scientist at Skolkovo Institute of Science and Technology.}
\cortext[cor1]{Corresponding author}
\author[add4,add2]{Alexander Novikov}
\ead{sasha.v.novikov@gmail.com}
\author[add3,add4]{Ivan Oseledets}
\ead{i.oseldets@skoltech.ru}

\address[add1]{Seminar for Applied Mathematics, ETH Zurich, R\"amistrasse 101, 8092 Zurich, Switzerland}
\address[add4]{Marchuk Institute of Numerical Mathematics of the Russian Academy of Sciences, 119333 Moscow, Russia}
\address[add2]{National Research University
Higher School of Economics, 101000 Moscow, Russia}
\address[add3]{Skolkovo Institute of Science and Technology, Skolkovo Innovation Center, 143026 Moscow, Russia}

\begin{abstract}
Such problems as computation of spectra of spin chains and vibrational spectra of molecules can be written as \emph{high-dimensional eigenvalue problems}, i.e., when the eigenvector can be naturally 
represented as a multidimensional tensor.
Tensor methods have proven to be an efficient tool 
for the approximation of solutions of high-dimensional eigenvalue problems, however, their performance deteriorates quickly when the number of eigenstates to be computed increases. We address this issue by designing a new algorithm motivated by the ideas of \emph{Riemannian optimization} (optimization on smooth manifolds)
for the approximation of multiple eigenstates in the \emph{tensor-train format}, which is also known as matrix product state representation.  
The proposed algorithm is implemented in TensorFlow, which allows for both CPU and GPU parallelization.

\end{abstract}

\end{frontmatter}

\section{Introduction} \label{sec:intro}

The paper aims at the approximate computation of $b$ lowest eigenvalues~$\energy^{(i)}$ and corresponding eigenvectors $\mX^{(i)}$ of a high-dimensional Hamiltonian $\mH\in\mathbb{R}^{n^d\times n^d}$:
\begin{equation}\label{eq:main}
	\mH \mX^{(i)} = \energy^{(i)}\mX^{(i)}, \quad i = 1,\dots, \bsize.
\end{equation}
Such problems arise in different applications of solid state physics and quantum chemistry problems including, but not restricted to the computation of spectra of spin chains ($n = 2$) or vibrational spectra of molecules (usually $n<20$).

High-dimensional problems are known to be computationally challenging due to the \emph{curse of dimensionality}, which implies exponential growth of the number of parameters with respect to the dimensionality~$d$ of the problem.
For example, storage of a single eigenvector for a spin chain with $d=50$ spins requires approximately 10 Pb of computer memory, which is far beyond available RAM on modern supercomputers.   
One way to deal with the curse of dimensionality is to utilize tensor decompositions of vectors $\mX^{(i)}$, $i=1,\dots,\bsize$ represented as multidimensional arrays $\tX^{(i)}$ of size $n\times \dots \times n$ also known as tensors.
Tensor decompositions have appeared to be useful for a long time both in solid states physics and molecular dynamics communities where different kinds of tensor decompositions have been used.
Recently, the tensor method has also attracted the attention in the numerical analysis community, where new ideas such as cross approximation method~\cite{ost-tucker-2008,ot-ttcross-2010,lars-htcross-2013} have been developed.

In this work we introduce a new method for solving problem~\eqref{eq:main} using the~\emph{tensor-train (TT) format}~\cite{osel-tt-2011}, also known as matrix product state (MPS) representation, which has been used for a long time in quantum information theory and solid state physics to
approximate certain wavefunctions \cite{white-dmrg-1992, ostlund-dmrg-1995} (DMRG method), see the review \cite{schollwock-2011} for more details.
One of the peculiarities of the proposed method is that it utilizes ideas of optimization on smooth manifolds (Riemannian optimization). 
This is possible due to the fact that the set of tensors of a fixed TT-rank\footnote{TT-rank controls the number of parameters in the decomposition (Sec.~\ref{sec:tt}).} 
 forms a smooth nonlinear manifold~\cite{holtz-manifolds-fixed-rank-2011}.
For small TT-ranks the manifold is low-dimensional and hence all computations are inexpensive.

To find a single eigenpair $(\energy^{(1)}, \mX^{(1)})$ Riemannian optimization allows to naturally avoid the rank growth of the method (Sec.~\ref{sec:riemtt}).
This is essential for fast computations as the complexity of the tensor method usually has strong rank dependence.
However, for finding several eigenvalues, the generalization is not straightforward and leads to a significant complexity increase.
Alternatively, one can compute eigenpairs using Riemannian optimization one-by-one, i.e. using deflation techniques, but such an approach is known to have slow convergence if the spectra is clustered~\cite{arbenz-lectures-2012}.
To avoid this we propose a modification of the method (Sec.~\ref{sec:proposed}) that on the one hand keeps the benefits of efficient single eigenstate computation using Riemannian approach, and on the other hand inherits faster convergence properties of block methods.

The idea of the proposed method is as follows. 
We suggest projecting all the eigenvectors at each iteration onto a single specifically chosen tangent space\footnote{In a nutshell, tangent space is a linearization of the manifold at a given point (Sec.~\ref{sec:riem_lopcg}).}. 
However, there is no reason to believe that all eigenvectors can be approximated using tangent space of a single eigenvector. 
Thus, in our algorithm, we always treat the projection as a correction to the current iterate.
Overall, this leads to a small, but non-standard optimization procedure for the coefficients of the iterative process (Sec.~\ref{sec:coef}). 

The main contributions of this paper are:
\begin{itemize}
	\item We develop a low-rank Riemannian alternating projection (LRRAP) concept for block iterative methods and apply it to the locally optimal block preconditioned conjugate gradient (LOBPCG) method. 
    \item We implement the proposed LRRAP LOBPCG solver using TensorFlow\footnote{If you are unfamiliar with TensorFlow, see~\ref{sec:app-t3f-implementation} for introduction.} library, which allows for parallelization on both CPUs and GPUs.
    \item We accurately calculate vibrational spectra of acetonitrile (CH$_3$CN) and ethylene oxide (C$_2$H$_4$O) molecules as well as spectra of one-dimensional spin chains. The comparison with the state-of-the-art methods in these domains indicates that we obtain comparable accuracy with a significant acceleration (up to 20 times) for large~$\bsize$  thanks to the design of the method and the GPU support.
\end{itemize}

\section{Riemannian optimization on TT manifolds} \label{sec:riemtt}

In this section we provide a brief description of TT-decomposition and Riemannian optimization essentials on the example of a single eigenvalue computation using LOBPCG method.
In Section~\ref{sec:proposed} the approach described here will be generalized to the computation of several eigenvalues.

\subsection{TT representation}\label{sec:tt}

Recall that we consider problem~\eqref{eq:main} and represent each eigenvector of size~$n^d$ as a $n\times \dots \times n$ tensor.
This allows for compression using tensor decompositions of multidimensional arrays, and particularly the TT decomposition. For tensor $\tX = \{\tX_{i_1,\dots,i_d}\}_{i_1,\dots,i_d=1}^{n}\in \mathbb{R}^{n\times \dots \times n}$ its TT decomposition reads
\begin{equation}\label{eq:tt-repr}
	\tX_{i_1,\dots,i_d} = G_{1}(i_1)\, G_{2}(i_2) \dots G_d (i_d),
\end{equation}
where $G_{k}(i_k)$ are $r_{k-1}\times r_k$ matrices, $k=2,\dots,d-1$.
For the product of matrices to be a number we require that $G_{1}(i_1)$ be row matrices $1\times r_1$ and $G_{d}(i_d)$ be column matrices $r_{d-1}\times 1$, which means $r_0=r_d = 1$.
For simplicity we force $r_1=\dots = r_{d-1} = r$ and call $r$ the \emph{TT-rank}. 
The same value of $r_k$ for all $k=1,\dots,d-1$ implies that for some modes $r_k$ can be overestimated.
Note that in order to store TT representation of array $\tX$ one needs $\mathcal{O}(dnr^2)$ elements compared with $n^d$ elements of the initial array.
The TT representation of Hamiltonian\footnote{
Matrix $\mH\in\mathbb{R}^{n^d\times n^d}$ can also be naturally considered as a multidimensional array $\tH$ of dimension $2d$.
The TT decomposition of $\tH$ reads 
$\tH_{i_1,\dots,i_d,j_1,\dots,j_d} = H_{1}(i_1,j_1)\, H_{2}(i_2,j_2) \dots H_d (i_d, j_d)$, where $i_1,\dots,i_d$ represent row indexing of $\mH$, while $j_1,\dots,j_d$ represent its column indexing, $H_{k}(i_k,j_k)$ are $R_{k-1}\times R_k$ matrices, $k=2,\dots,d-1$ and $R_0=R_d = 1$.
} 
is defined by analogy.
We will denote the maximum rank of a Hamiltonian by $R$. 

It rarely happens that some tensor can be represented with small TT-rank~$r$ exactly.
Therefore, to keep ranks small the tensor is approximated with some accuracy by another tensor with a small TT-rank.
In the considered numerical experiments we expect exponential decay of the introduced error with respect to~$r$.

From now on we say that vector $\mX$ of length $n^d$ is of TT-rank $r$ implying that being reshaped into a $n\times\dots\times n$ multidimensional array $\tX$: $\mX=\mathrm{vec}(\tX)$ it can be represented with TT-rank equal to $r$.

\subsection{Rayleigh quotient minimization using Riemannian optimization}\label{sec:riem_lopcg}

Consider the problem of finding the smallest eigenvalue~$\energy^{(1)}$ (assuming it is simple) and the corresponding eigenvector~$\mX^{(1)}$.
Suppose we are given an a priori knowledge that $\mX^{(1)}$ can be accurately approximated by a TT representation of a small TT-rank $r$, i.e. it lies on the manifold of tensors of fixed TT-rank $r$:
\[
	\M_r \equiv \{\mX\in \Rnd \, |\, \text{TT-rank}(\mX) = r \}.
\]
To obtain the eigenpair $(\energy^{(1)},\mX^{(1)})$ we pose a Rayleigh quotient $\RQ(\mX)$ minimization problem on the low-parametric manifold $\M_r$ instead of the full $\Rnd$ space:
\begin{equation}\label{tfeig:eq:rqlr}
	\min_{\mX \in \M_r}\, \RQ(\mX), \quad  \RQ(\mX)\equiv \frac{\left< \mX, \mH \mX \right>}{\left< \mX, \mX \right>},
\end{equation}
where we assume that $\mH$ is symmetric.
Since $\M_r$ forms a smooth manifold~\cite{holtz-manifolds-fixed-rank-2011} 
one can utilize the so-called methods of Riemannian optimization, i.e. optimization on smooth manifolds.

One of the key concepts for the Riemannian optimization is \emph{tangent space}, which consists of all tangent vectors to $\M_r$ at a given~$\mX$~\cite{robsal-diffgeom-2011}.
We will denote tangent space of $\M_r$ at $\mX$ by $T_\mX \M_r$.
tangent space can be viewed as the linearization of a manifold at the given point $\mX$. 
It has the same dimension as the manifold~\cite{robsal-diffgeom-2011} and assuming that $r$ is small, it allows to locally replace the manifold with a low-dimensional linear space.

\begin{figure}[t]
\begin{center}
\begin{tikzpicture}
\node at (0,0) {\includegraphics[width=100mm]{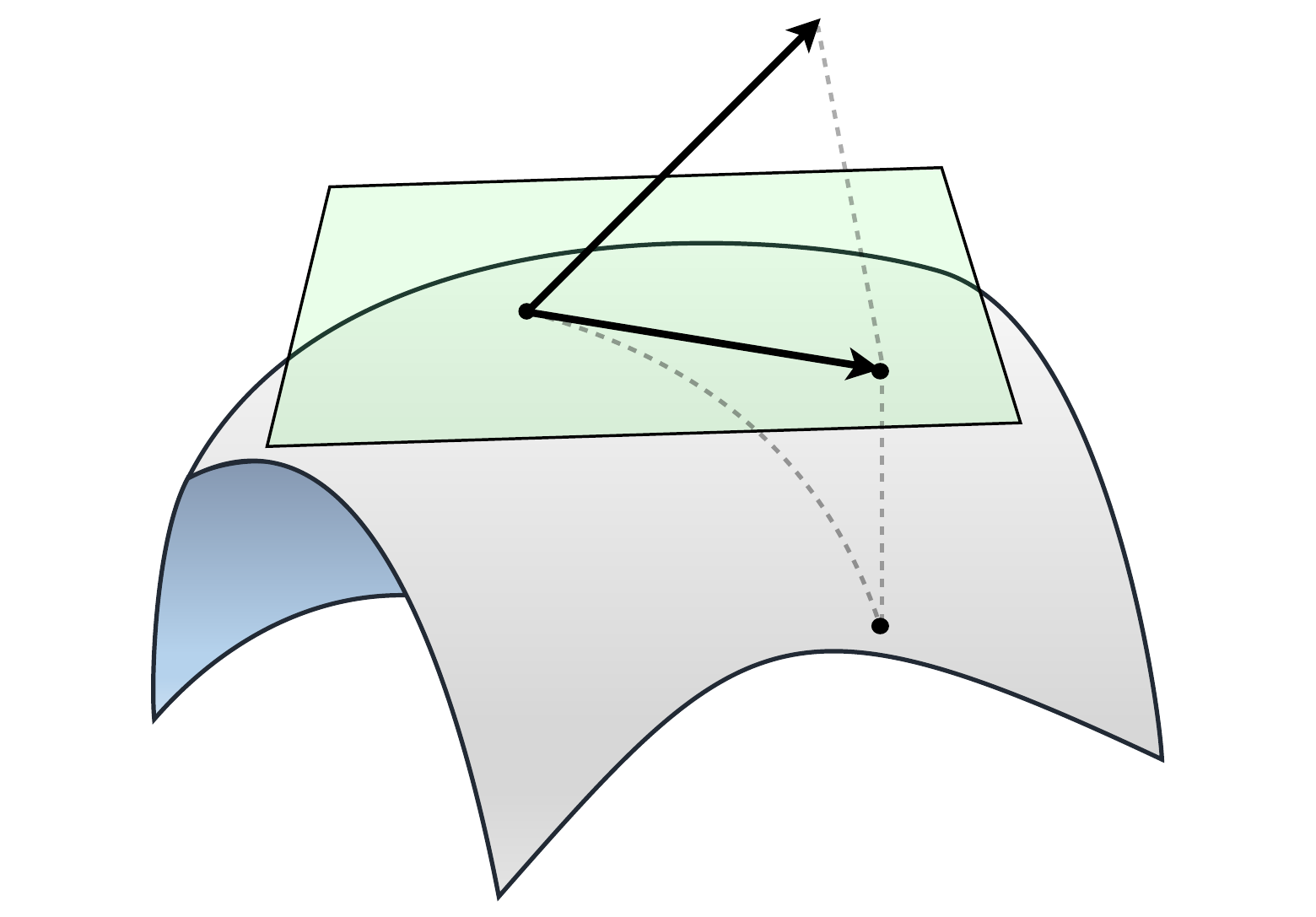}};
\node at (-1,-2.2) {{ $\M_r$}};
\node at (-1.9,1.8) {{ $T_{\mX_k} \M_r$}};
\node at (-0.15,3.2) {{ $-\nabla \RQ (\mX_k)$}};
\node [rotate=-7] at (0.5,1.17) {{\small $-\OP_{\mX_k}\negthickspace \nabla \RQ (\mX_k)$}};
\node at (2.2,-1.3) {{ $\mX_{k+1}$}};
\node at (-1.3,0.9) {{ $\mX_{k}$}};
\end{tikzpicture}
\end{center}
\caption{Illustration of the Riemannian gradient descent method. $\M_r$ denotes the smooth manifold of vectors of fixed TT-rank and  $T_{\mX_k} \M_r$ its tangent space at $\mX_{k}$. The gradient is projected on the tangent space and then $\mX_k$ is moved in the direction of the projected gradient and afterwards is retracted to the manifold: $\mX_{k+1} = \Retr_r \left(\mX_k - \OP_{\mX_k} \nabla\RQ (\mX_k) \right)$. 
} \label{fig:riem}
\end{figure}

Provided a tangent space at hand we can discuss the simplest optimization method on $\M_r$ --- the Riemannian gradient descent method, which consists of several steps and is illustrated in Fig.~\ref{fig:riem}.
	Given the starting point $\mX_k$, the first step is to calculate the Riemannian gradient, which in this case is projection of the standard Euclidean gradient: $\OP_{\mX_k} \nabla \RQ (\mX_{k})$, where $\OP_{\mX_k}$ denotes the orthogonal projection on $T_{\mX_k} \M_r$. 
	This is simply the Rayleigh quotient steepest descent direction at $\mX_k$, restricted to the tangent space~$T_{\mX_k} \M_r$.
The second step is, given the search direction from $T_\mX \M_r$, map it back to the manifold to obtain $\mX_{k+1}$ of the iterative process.
This is done by the smooth mapping $\widetilde\Retr_r(\mX,\cdot):T_\mX \M_r\to \M_r$ called \emph{retraction}.
Note that in case of fixed-rank manifold $\M_r$ we use the retraction satisfying $\widetilde\Retr_r(\mX,\xi) \equiv \Retr_r(\mX + \xi)$~\cite{ao-retract-2014} and call it truncation.

Similarly to the original gradient descent method, the convergence of such method can be slow and preconditioning has to be used.
There are different ways how to account for a preconditioner in the Riemannian version of the iteration.
We will use the idea considered in \cite{ksv-manprec-2016}, where the preconditioner acts on the gradient and the result is projected afterwards
\begin{equation}\label{eq:riem_pinvit}
\begin{split}
	 \mX_{k+1} =  \Retr_r \left(\mX_{k} - \,\OP_{\mX_k} \Prec \nabla\RQ (\mX_k) \right), \\ 
\end{split}
\end{equation}
which for eigenvalue computations can be viewed as Riemannian generalization\footnote{The search direction~$\OP_{\mX_k} \Prec \nabla\RQ (\mX_k)$ is usually additionally multiplied by a constant $\rel_k$ to be found from the line search procedure $\RQ(\mX_k-\rel_k \OP_{\mX_k} \Prec \nabla\RQ (\mX_k))\to\min_{\rel_k}$. This ensures the convergence in the presence of the nonlinear mapping $\Retr_r$.} of the preconditioned inverse iteration (PINVIT).
Indeed, one iteration of a classical version of PINVIT reads
\[
	 \mX_{k+1} = \mX_{k} - \Prec \Rvec_k, 
\]
where $\Rvec_k = \mH \mX_k - \RQ (\mX_k) \mX_k$ denotes the residual, which is proportional to the gradient of the Rayleigh quotient:
\[
	\nabla \RQ (\mX) = \frac{2}{\left<\mX, \mX \right>} \left(\mH \mX - \RQ (\mX) \mX \right).
\]
Note that to avoid growth of $\|\mX_{k}\|$ additional normalization $\mX_{k}:=\mX_{k}/\|\mX_{k}\|$ is done after each iteration to ensure $\left<\mX_k, \mX_k \right>=1$.

We could have restricted ourselves to the case of PINVIT~\eqref{eq:riem_pinvit}, but to get faster convergence we utilize a superior method --- locally optimal preconditioned conjugate gradients (LOPCG) and its block version (LOBPCG) to calculate several eigenvalues (see Sec.~\ref{sec:proposed}).
To our knowledge the Riemannian version of LOPCG was not considered in the literature, so we  provide it here.
According to the classical LOBPCG the search direction $\Pvec_k$ is a linear combination of the preconditioned gradient and $\Pvec_{k-1}$.
In the Riemannian setting instead of the preconditioned gradient $\Prec \nabla\RQ (\mX_k)$ we consider its projected analog~$\OP_{\mX_k} \Prec \nabla\RQ (\mX_k)\in T_{\mX_k} \M_r$.
However, the problem is that~$\Pvec_{k-1}\not\in T_{\mX_k} \M_r$, so similarly to~\cite{v-matcompl-2013} we use another important concept from the differential geometry -- vector transport, which is a mapping from $T_{\mX_{k-1}}\M_r$ to $T_{\mX_{k}}\M_r$ satisfying certain properties~\cite{absil}. 
Since $\M_r$ is an embedded submanifold of $\mathbb{R}^{n^d}$, the orthogonal projection from~$T_{\mX_{k-1}}\M_r$ to $T_{\mX_{k}}\M_r$ can be used as a vector transport~\cite{absil}, so that
\begin{equation}\label{eq:pk1}
	\Pvec_{k+1} = c_1 \,\OP_{\mX_k} \Prec \nabla \RQ (\mX_k) + c_2 \,\OP_{\mX_k} \Pvec_{k},
\end{equation}
and 
\begin{equation}\label{eq:lopcg_riem}
	\mX_{k+1} = \Retr_r \left( \mX_k + \Pvec_{k+1} \right)
\end{equation}
where $c_1, c_2$ are constants to be found from 
\begin{equation}\label{eq:c1c2}
	(c_1,c_2) = \argmin_{(\zeta_1,\zeta_2)} \RQ (\mX_k + \zeta_1 \,\OP_{\mX_k} \Prec \nabla \RQ (\mX_k) + \zeta_2 \,\OP_{\mX_k} \Pvec_{k})
\end{equation}
Note that we have omitted $\Retr_r$ in the optimization procedure.
This allows to solve~\eqref{eq:c1c2} exactly.
Despite in numerical computations we have never faced the problem of functional increase with the omitted $\Retr_r$, additional line-search procedure could be introduced to ensure convergence.

\subsection[Complexity reduction for computations on Riemannian manifold]{Complexity reduction for computations on $\M_r$}

One of the key benefits of the Riemannian optimization on $\M_r$ is that it allows to avoid the rank growth.
This is particularly important in low-rank tensor computations due to the strong rank dependence of tensor methods, which sometimes is called the curse of the rank.
The crucial property that allows for a significant speed-up of computations in the Riemannian approach is that TT-ranks of any vector from a tangent space of $\M_r$ has ranks at most $2r$.
Since tangent space is a linear space, linear combination of any number of vectors from one tangent space also has rank at most $2r$.
As a result, $\Pvec_{k+1}$ in \eqref{eq:pk1} is of rank at most $2r$ and for the already computed $\Pvec_{k+1}$ the computation of~\eqref{eq:lopcg_riem} is inexpensive.
By contrast, if in~\eqref{eq:lopcg_riem} we omit $\OP_{\mX_k}$, i.e. no Riemannian optimization approach is used, then this leads to complexity increase since matrix-vector multiplications considerably increase TT-ranks.
Moreover, one can show that in~\eqref{eq:c1c2} to find $c_1$ and $c_2$ scalar products of TT-tensors have to be computed. The calculation of scalar products of two vectors belonging to the same tangent space is of less complexity than of two general TT-tensors of the same rank.
We provide details about the implementation of the Riemannian optimization for TT-manifolds in Sec.~\ref{sec:impl} after general description of the proposed method.

\section{The proposed method} \label{sec:proposed}

The main goal of the paper is to consider a problem of finding several eigenvalues and the corresponding eigenvectors.
For convenience we rewrite~\eqref{eq:main} in the block form
\[
	\mH \X = \X\Lam, \quad \X = [\mX^{(1)},\dots,\mX^{(\bsize)}], \quad \Lam = \mathrm{diag}(\energy^{(1)},\dots,\energy^{(\bsize)}),
\]
additionally assuming $\energy^{(1)}\leq \energy^{(2)} \leq \dots \leq \energy^{(\bsize)}\not=\energy^{(\bsize+1)}$.
Then the problem under the TT-rank constraint on $\mX^{(\alpha)}$, $\alpha=1,\dots,\bsize$ can be reformulated as trace minimization:
\begin{equation}
\label{eq:lobpcg_riem_problem_formulation}
\begin{aligned}
& \underset{\multipleX\in \mathbb{R}^{n^d\times \bsize}}{\text{minimize}}
& & \Tr(\multipleX^\intercal \mH \multipleX) \\
& \text{subject to}
& & \multipleX^\intercal \multipleX = \I_\bsize\\
& & & \mX^{(i)} \in \M_\rank, \; i = 1, \ldots, \bsize,
\end{aligned}
\end{equation}
where $\Tr(\cdot)$ denotes trace of a matrix.

\subsection{Speeding-up computation of the smallest eigenvalue}

As a starting point to generalize~\eqref{eq:lopcg_riem} to the block case assume that we are first aiming at finding the first eigenvalue.
The iteration~\eqref{eq:lopcg_riem} can be improved by performing additional subspace acceleration by using more vectors in the tangent space of $\M_r$ at $\mX_{k}^{(1)}$.
In particular, we suggest using the LOBPCG method and project all arising vectors to the tangent space\footnote{Note that $\OP_{\mX_k^{(1)}} \mX_k^{(1)} = \mX_k^{(1)}$.} at $\mX_k^{(1)}$:
\begin{equation}
\label{eq:p1_lobpcg_riem}
\begin{split}
	& \R_k = \mH \X_k - \X_k\Lam_k, \quad \Lam_k = \mathrm{diag}\left(\RQ\left(\mX_k^{(1)}\right),\dots, \RQ\left(\mX_k^{(\bsize)}\right)\right) \\
    & \PP_{k+1} = \OP_{\mX^{(1)}_k} \Prec \R_k \C_2  + \OP_{\mX^{(1)}_k} \PP_{k}\C_3\\
    & \X_{k+1} =\Retr_r \left(\OP_{\mX_k^{(1)}} \X_k \C_1 + \PP_{k+1} \right),
\end{split}
\end{equation}
where the truncation operator $\Retr_r$ is applied independently to each column of the matrix and $\C_i\in\mathbb{R}^{\bsize \times \bsize}$, $i=1,2,3$ are matrices of coefficients to be found from trace minimization.
In particular, introducing notation 
\[
	\basis_k = \OP_{\mX_k^{(1)}}\begin{bmatrix}\X_{k} &\R_k &\PP_{k}\end{bmatrix}\in \mathbb{R}^{n^d\times 3\bsize}, \quad \C = \begin{bmatrix} \C_1 \\ \C_2 \\ \C_3 \end{bmatrix} \in \mathbb{R}^{3 \bsize \times \bsize},
\]
we have
\begin{equation}
\label{eq:naive_coef_problem_lobpcg}
\begin{aligned}
& \underset{\C}{\text{minimize}}
& & \Tr(\C^\intercal(\basis_k^\intercal \mH \basis_k)\C), \\
& \text{subject to}
& & \C^\intercal(\basis_k^\intercal\basis_k)\C = \I_\bsize,
\end{aligned}
\end{equation}
which reduces to a classical generalized eigenvalue problem of finding $\bsize$ smallest eigenvalues $\theta_1,\dots,\theta_\bsize$ and corresponding eigenvectors of the matrix pencil $\basis_k^\intercal \mH \basis_k - \theta \basis_k^\intercal\basis_k$ of size $3\bsize\times 3\bsize$.

Since we use projection to the same tangent space, there is no rank growth even for large $\bsize$ and hence the application of~$\Retr_r$ is inexpensive.
Moreover, we expect that to achieve a given accuracy of $\energy^{(1)}$ the number of iterations for~\eqref{eq:p1_lobpcg_riem} is less than the number of iterations for~\eqref{eq:lopcg_riem} thanks to the additional subspace acceleration.

\subsection{Riemannian alternating projection method}

Similarly to~\eqref{eq:p1_lobpcg_riem} all column vectors of $[\X_{k},\Prec \R_k,\PP_{k}]$ can be projected to the tangent space $T_{\mX^{(t_k)}_k}\M_r$ of $\mX^{(t_k)}_k$ for some integer $1\leq t_k\leq \bsize$, which not necessarily equals $1$.
However, there is no evidence that all eigenvectors except for the $t_k$-th one can be accurately approximated using $T_{\mX^{(t_k)}_k}\M_r$.
Therefore, we propose to use $T_{\mX^{(t_k)}_k}\M_r$ to search for the \emph{correction} to the already found approximations of $\X_k$ as follows:
\begin{equation}
\label{eq:lobpcg_riem}
\begin{split}
	& \R_k = \mH \X_k - \X_k\Lam_k, \\
    & \PP_{k+1} = \OP_{\mX^{(t_k)}_k} \Prec \R_k \C_2  + \OP_{\mX^{(t_k)}_k} \PP_{k}\C_3\\
    & \X_{k+1} =\Retr_r \left( \textcolor{blue}{\X_k\,\diag(\brel)} + \OP_{\mX^{(t_k)}_k} \X_k \C_1 + \PP_{k+1} \right),
\end{split}
\end{equation}
where $\brel\in\mathbb{R}^{\bsize}$ an $\diag(\brel)$ denotes diagonal $\bsize\times \bsize$ matrix with the vector $\brel$ on the diagonal.
By forcing $\X_k$ to be multiplied by a diagonal matrix $\diag(\brel)$ instead of a general square matrix, we reduce computational cost of the method.
Indeed, in this case we avoid calculating linear combinations of the columns of $\X_k$ which do not belong to the same tangent space.

Note that we may also write
\[
	\X_{k+1} =\Retr_r \left( \X_k\, \diag(\brel) +  \basis \C \right),
\]
where
\[
	\basis_k = \OP_{\mX^{(t_k)}_k}\left[ \X_{k}, \Prec \R_k, \PP_{k}\right],
\]
and the matrices of coefficients
$$
\C = \begin{bmatrix} \C_1 \\ \C_2 \\ \C_3 \end{bmatrix},\quad \C_j\in\mathbb{R}^{\bsize\times \bsize}, \quad j=1,2,3 
$$ 
and $\brel\in\mathbb{R}^{\bsize}$ are to be found from the trace minimization problem
\begin{equation}
\label{eq:coefficients-problem0}
\begin{aligned}
& \underset{\brel,\, \C}{\text{minimize}}
& & \Tr\left[(\X_k\, \diag(\brel) + \basis_k \C)^\intercal \mH (\X_k\, \diag(\brel) + \basis_k \C) \right]\\
& \text{subject to}
& & (\X_k\, \diag(\brel) + \basis_k \C)^\intercal (\X_k \, \diag(\brel) + \basis_k \C) = \I_{\bsize}
\end{aligned}
\end{equation}
or equivalently
\begin{equation}
\label{eq:coefficients-problem}
\begin{aligned}
& \underset{\brel,\, \C}{\text{minimize}}
& & \Tr\left( 
\begin{bmatrix}
\mathrm{diag}(\brel) \\
\C	
\end{bmatrix}^\intercal
\begin{bmatrix}
\X_k^\intercal \mH \X_k & \X_k^\intercal \mH \basis_k \\
\basis_k^\intercal \mH \X_k & \basis_k^\intercal \mH \basis_k
\end{bmatrix}
\begin{bmatrix}
\mathrm{diag}(\brel) \\
\C	
\end{bmatrix}
\right)\\
& \text{subject to}
& & \begin{bmatrix}
\mathrm{diag}(\brel) \\
\C	
\end{bmatrix}^\intercal
\begin{bmatrix}
\X_k^\intercal  \X_k & \X_k^\intercal  \basis_k \\
\basis_k^\intercal  \X_k & \basis_k^\intercal  \basis_k
\end{bmatrix}
\begin{bmatrix}
\mathrm{diag}(\brel) \\
\C	
\end{bmatrix} = \I_{\bsize}
\end{aligned}
\end{equation}
which because of the diagonal constraint does not boil down to a standard generalized eigenvalue problem.
Due to the description technicality of solution of~\eqref{eq:coefficients-problem}, we postpone it to Section~\ref{sec:coef}.

The convergence of the proposed method depends on the choice of the integer sequence $\{t_1,t_2,\dots, t_n, \dots\}$. 
We call the strategy that in a certain way chooses the tangent space on each iteration the \emph{tangent space schedule}. 
Note that we could have found all eigenvalues one by one, i.e. $t_1=\dots = t_{k_1^{\delta}}=1$, then $t_{k_1^{\delta}}=\dots = t_{k_1^{\delta} + k_2^{\delta}}=2$ and so on, where $k_\alpha^{\delta}$ is the number of iterations for $\RQ(\mX^{(\alpha)})$ to achieve accuracy $\delta$.
This strategy, however, for a large number of eigenvalues~$\bsize$ requires a lot of iterations to be done.
Therefore, we utilize strategies that do not require all tangent spaces to be chosen at least ones.
We found that although the random choice (discrete uniform distribution) of $t_i$: $1\leq t_i \leq \bsize$ already ensures convergence in most of the cases, the strategy
\begin{equation}
t_k = \argmax_{i} \abs{\frac{ \RQ (\mX_{k-1}^{(i)}) - \RQ (\mX_k^{(i)})}{\RQ (\mX_k^{(i)})}}
\label{eq:schedule}
\end{equation}
in which we choose eigenvalue with the current slowest convergence, yields more reliable results for larger number of eigenvalues.
Before running the adaptive strategy for $t_k$ we choose $t_1=\dots = t_{k_0}=1$ while convergence criteria for $\RQ(\mX^{(1)})$ is not fulfilled.
This corresponds to~\eqref{eq:p1_lobpcg_riem} instead of~\eqref{eq:lobpcg_riem} and hence, speeds up computations.
The implementation details of the iteration~\eqref{eq:lobpcg_riem} will be given in the next section and the algorithm for the computation of coefficients will be described later in Section~\ref{sec:coef}.

\section{Implementation of tensor operations} \label{sec:impl}

In this section we provide a brief description of implementation details for Riemannian optimization on $\M_r$ with complexity estimates and summarize the algorithm. 

\subsection{Representation of tangent space vectors}
We start with the description of vectors from tangent spaces, as they are the main object we are working with.
Let $\mX$ be given by its TT decomposition~\eqref{eq:tt-repr}.
It is known~\cite{osel-tt-2011} that by a sequence of QR decomposition of cores, \eqref{eq:tt-repr} can be represented both as
\begin{equation}\label{eq:uqr}
	\tX_{i_1\dots i_d} = U_1(i_1) U_2 (i_2) U_3 (i_3) \dots U_d (i_d),
\end{equation}
where $\mathcal{U}_{k}=\{U_k(i_k)\}_{i_k=1}^n\in\mathbb{R}^{r_{k-1}\times n\times r_{k}}$, $k=1,\dots,d-1$ being reshaped into matrices $\mathscr{M}^{\mathsf{L}}_k(\mathcal{V}_{k})$ of size $r_{k-1} n\times r_{k}$ have orthogonal columns, $r_1=\dots = r_{d-1}=r$ and $r_0=r_d=1$.
Here $\mathscr{M}_k^{\mathsf{L}}:\mathbb{R}^{r_{k-1}\times n\times r_{k}}\to \mathbb{R}^{r_{k-1} n\times r_{k}}$ denotes the matricization operator, that maps first two indices of considered three-dimensional arrays into a single one.
Similarly, we may have
\begin{equation}\label{eq:vqr}
	\tX_{i_1\dots i_d} = V_1(i_1) V_2 (i_2) V_3 (i_3) \dots V_d (i_d),
\end{equation}
where $\mathcal{V}_{k}=\{V_k(i_k)\}_{i_k=1}^n\in\mathbb{R}^{r_{k-1}\times n\times r_{k}}$, $k=2,\dots,d$ being reshaped into matrices $\mathscr{M}^{\mathsf{R}}(\mathcal{V}_{k})$ of size $r_{k-1} \times nr_{k}$ have orthogonal rows.
Representations~\eqref{eq:uqr} and~\eqref{eq:vqr} are called correspondingly left- and right-orthogonalizations of TT-representation~\eqref{eq:tt-repr}, which can be performed with $\mathcal{O}(dnr^3)$ complexity.
Using these notations, one possible way to parametrise tangent space $T_{\mX}\M_r$ is as follows.
Any $\xi = \mathrm{vec}(\Xi)\in T_\mX \M_r$ can be represented as
\begin{equation}\label{eq:tangent_param}
\begin{split}
	\Xi_{i_1\dots i_d} = 
				&\delta G_1(i_1) V_2 (i_2) V_3 (i_3) \dots V_d (i_d)\ + \\
				&U_1(i_1) \delta G_2 (i_2) V_3 (i_3) \dots V_d (i_d)\ + \dots + \\
				&U_1(i_1) U_2 (i_2) U_3 (i_3) \dots \delta G_d (i_d),
\end{split}
\end{equation}
where cores $\delta\mathcal{G}_{k}=\{\delta G_k(i_k)\}_{i_k=1}^n\in\mathbb{R}^{r_{k-1}\times n\times r_{k}}$ additionally satisfy the following gauge conditions to ensure uniqueness\footnote{To see that the tangent space is overparametrized by~\eqref{eq:tangent_param} note that the number of parameters in all the $G_k(i_k)$ is $(d-2)nr^2 + 2nr$, while the dimension of the manifold~$\M_r$ and hence of the tangent space is smaller: $\mathrm{dim}(T_\mX \M_r)=(d-2)nr^2 + 2nr - (d-1)r^2$~\cite{holtz-manifolds-fixed-rank-2011}.} of the representation:
\begin{equation}\label{eq:gauge}
	\left(\mathscr{M}^{\mathsf{L}}\left(\delta\mathcal{G}_{k} \right)\right)^\intercal \mathscr{M}^{\mathsf{L}}\left(\mathcal{U}_{k} \right) = \mathbf{0}, \quad k=1,\dots,d-1.
\end{equation}
To show that any vector from a tangent space has TT-rank not greater than $2r$ we note that for~\eqref{eq:tangent_param} the explicit formula holds:
\[
	\Xi_{i_1\dots i_d} = S_1(i_1) S_2(i_2) \dots S_d (i_d),
\]
where for $k=2,\dots,d-1$
\[
S_1(i_1) = 
	\begin{bmatrix}
			\delta G_1(i_1) & U_1(i_1)
	\end{bmatrix},
	\,
S_k(i_k) = 
	\begin{bmatrix}
			V_k(i_k) &  \\
			\delta G_k(i_k) & U_k (i_k)
	\end{bmatrix},
	\,
S_d(i_d) = 
	\begin{bmatrix}
			 V_d(i_d) \\ \delta G_d(i_d)
	\end{bmatrix},
\]
which can be verified by direct multiplication of all $S_k (i_k)$.
From~\eqref{eq:tangent_param} it is also simply noticeable that if $\xi^{(1)}, \xi^{(2)}\in T_\mX \M_r$ given by $\delta G^{(1)}_k(i_k)$ and $\delta G^{(2)}_k(i_k)$, $k=1,\dots,d$ correspondingly, then their linear combination $(\alpha\xi^{(1)} + \beta\xi^{(2)})\in T_\mX \M_r$ is given by $(\alpha\, \delta G^{(1)}_k(i_k) + \beta\, \delta G^{(2)}_k(i_k))$, $\alpha,\beta\in\mathbb{R}$.

\subsection{Projection to a tangent space}\label{sec:projection}

Representation~\eqref{eq:tangent_param} allows to obtain explicit formulas for $\delta G_k(i_k)$ of $\xi = \OP_{T_\mX \M_r}\mZ$ without inversions of possibly ill-conditioned matrices~\cite{ksv-manprec-2016}: 
\[
\begin{split}
	&\mathrm{vec}(\delta\mathcal{G}_k)= \left(\I_{r_k-1}\otimes (\I_{nr_k} - 
	\mathscr{M}^{\mathsf{L}}\left(\mathcal{U}_{k} \right)
	\mathscr{M}^{\mathsf{L}}\left(\mathcal{U}_{k} \right)^\intercal)\right)
	(\mathbf{X}_{>k}^\intercal \otimes I_n \otimes \mathbf{X}_{<k}^\intercal)\, \mZ,\\
	&\mathrm{vec}(\delta\mathcal{G}_d)= 
	(I_n \otimes \mathbf{X}_{<d}^\intercal)\, \mZ,
\end{split}
\]
where $k=1,\dots, d-1$ and
\[
\begin{split}
	&\mathbf{X}_{<k} = [G_1(i_1)\dots G_{k-1}(i_{k-1})]\in \mathbb{R}^{n^{k-1} \times r}, \\
	&\mathbf{X}_{>k} = [G_{k+1}(i_{k+1})\dots G_{d}(i_d)]\in \mathbb{R}^{n^{d-k} \times r}.
\end{split}
\]
Matrices $\mathbf{X}_{<k}$, $\mathbf{X}_{>k}$ are never formed explicitly and used here for the ease of notation.
The complexity of projecting a vector~$\mZ$ given by its TT decomposition with $\text{TT-rank}(\mZ)= r_\mZ$ is $\mathcal{O}(dnr r_\mZ^2)$.

One of the most time-consuming operations arising in the algorithm is a projection to a tangent space of a matrix-vector product.
Suppose that both $\mH$ and $\mY$ are given in the TT format with TT-ranks $R$ and $r_\mY$.
Then $\mZ = \mH\mY$ can also be represented in the TT format with the TT-rank bounded from above as $r_\mY R$~\cite{osel-tt-2011} since
\[
\begin{split}
	&\tZ_{i_1,\dots,i_d} = Z_1(i_1)\dots Z_d(i_d), \\ 
	&Z_k(i_k) = \sum_{j_k=1}^n H_k(i_k, j_k) \otimes G_k(i_k) \in \mathbb{R}^{R_{k-1}(r_\mY)_{k-1}\times R_{k}(r_\mY)_{k}}.
\end{split}
\]
Calculation of $Z_k(i_k)$ can be done in $\mathcal{O}(n^2 R^2 r_\mY^2)$ complexity. Once $Z_k(i_k)$ $k=1,\dots,d$ are calculated, complexity of finding the TT representation of $\OP_{T_\mX \M_r} \mZ$ is $\mathcal{O}(dnr r_\mZ^2) = \mathcal{O}(dnr r_\mY^2 R^2)$ as $r_\mZ= r_\mY R$.
This complexity can be additionally reduced down to $\mathcal{O}(dn^2 r r_\mY R^2)$  if we take care of the Kronecker-product structure of $Z_k(i_k)$.
Moreover, since it does not require computation of $Z_k(i_k)$ explicitly, the storage is also less in this case.

\subsection{Computation of inner products} \label{sec:inner}
 
Next thing that arises in the method, in particular in~\eqref{eq:coefficients-problem}  is the computation of Gram matrices $\basis^\intercal \basis$, $\basis^\intercal \mH \basis$, $\bX^\intercal \bX$, $\basis^\intercal \bX$, $\basis^\intercal \mH \basis$, $\basis^\intercal \mH \X$  and as we will see in Sec.~\ref{sec:coef} solution to~\eqref{eq:coefficients-problem} will involve only $\mathrm{diag}(\bX^\intercal\mH \bX)$ instead of the full $\bX^\intercal\mH \bX$.
Here $\basis$ consists of columns that belong to a single tangent space $T_\mX \M_r$ and $\bX$ are general tensors of TT-rank $r$.

Let us first discuss the computation of $\basis^\intercal \basis$.
Let $\vv^{(\alpha)}$ be the $\alpha$-th column of~$\basis$, $\alpha=1,\dots,3\bsize$ and be given by $\delta G_k (i_k) \equiv \delta G_k^{(\alpha)} (i_k)$ from~\eqref{eq:tangent_param}. 
Then, thanks to gauge conditions~\eqref{eq:gauge}, left orthogonality of $\mathcal{U}$ and right orthogonality of $\mathcal{V}$, we have
\begin{equation}\label{eq:inner_tangent}
	\left< \vv^{(\alpha)}, \vv^{(\beta)} \right> = \sum_{k=1}^d  \left<\delta\mathcal{G}_k^{(\alpha)},  \delta\mathcal{G}_k^{(\beta)} \right>_F,
\end{equation}
where $\left<\mathcal{A},  \mathcal{B} \right>_F \equiv \mathrm{vec}(\mathcal{A})^\intercal \, \mathrm{vec}(\mathcal{B})$ is the Frobenius inner product.
Therefore, the complexity of computing $\basis^\intercal \basis$ is $\mathcal{O}(\bsize^2 dnr^2)$.
The computation of~$\basis^\intercal \mH \basis$ is done using the following trick.
Note that all columns of $\basis$ belong to a single tangent space $T_\mX \M_r$.
Hence, 
\begin{equation}\label{eq:trick_proj}
	\basis^\intercal \mH \basis = (\OP_{T_\mX \M_r}\basis)^\intercal \mH\basis = \basis^\intercal (\OP_{T_\mX \M_r}\mH \basis).
\end{equation}
As a result, we first compute $\OP_{T_\mX \M_r}\mH \basis$ using the procedure described in Section~\ref{sec:projection}.
Then we compute the inner product of two vectors from the same tangent space as in~\eqref{eq:inner_tangent}.
Thus, the complexity of finding $\basis^\intercal \mH \basis$ is $\mathcal{O}(\bsize dn^2 r^2 R^2 + \bsize^2 d r^2 )$.
Similarly to~\eqref{eq:trick_proj} we have $\basis^\intercal \mH \X = \basis^\intercal (\OP_{T_\mX \M_r}\mH \X)$ and $\basis^\intercal \X = \basis^\intercal (\OP_{T_\mX \M_r} \X)$.

Computation of $\bX^\intercal \bX$ is a standard procedure, which consists of $\bsize^2$ inner products of TT tensors. 
Since the complexity of inner product of two tensors of TT-rank $r$ is $\mathcal{O}(dnr^3)$~\cite{osel-tt-2011}, computing $\bsize^2$ of inner products require~$\mathcal{O}(\bsize^2 dnr^3)$ floating point operations.

Finally, the computation of ${\mX^{(\alpha)}}^\intercal \mH {\mX^{(\alpha)}}$, $\alpha=1,\dots,\bsize$ arising in $\mathrm{diag}(\bX^\intercal\mH \bX)$ costs $\mathcal{O}(nd R r^2 (r+R))$. 
Eventually, computation of $\mathrm{diag}(\bX^\intercal\mH \bX)$ has complexity $\mathcal{O}(\bsize nd R r^2 (r+nR))$.

\subsection{Retraction computation}
The final thing remains to compute to proceed to the next iteration after we have found the coefficients $\brel, \C$ and $\basis_k$, is to retract the obtained vectors $\X_{k}\mathrm{diag}(\brel) + \basis_k \C$ to the manifold $\M_r$.
All columns $\basis_k$ are from the same tangent space at $\mX_k^{(t_k)}$, so the correction $(\basis_k \C)[:,i]$ to each $\mX_k^{(i)}$, $i=1,\dots,\bsize$ is of rank $2r$.
Therefore, $\mX_k^{(i)} \brel[i] + (\basis_k \C)[:,i]$ is of rank at most $2r$ for $i= t_k$ and $3r$ otherwise. 
The retraction is done by the TT-SVD algorithm~\cite{osel-tt-2011} that consists of a sequence of QR and SVD decompositions applied to the unfolded tensor cores.
It has complexity $\mathcal{O}(dnr^3)$.

Note that the retraction is properly defined only for $i= t_k$, but we formally apply it to other vectors as well, as the coefficients $\brel,\C$ were obtained to ensure the descent direction for all vectors.
The descent direction is, however, chosen without regard to the retraction.
This can be accounted for by introducing approximate line search, but numerical experiments showed this is in general redundant and $\brel,\C$ already provide good enough approximation.

\subsection{The algorithm description}

Let us now discuss the iteration~\eqref{eq:lobpcg_riem} step-by-step.
For simplicity let us denote $\OP_k \equiv \OP_{\mX^{(t_k)}_k}$.
First, assuming we are given $\PP_{k}$, the calculation of $\OP_k \PP_{k}$ is done in $\mathcal{O}(\bsize dnr^3)$ complexity.

To calculate $\OP_k\Prec \R_k$ term we split it into two parts:
\[
	\OP_k\Prec \R_k = \OP_k\Prec\mH \X_k -\OP_k \Prec\X_k\Lam_k.
\]
The $\OP_k \Prec\X_k$ is a projected matrix vector product, which is calculated as described in Section~\ref{sec:projection} and costs $\mathcal{O}(\bsize dn^2 r^2 R^2)$ operations.
The term $\OP_k\Prec\mH \X_k$ is more difficult to compute as it involves two sequential matrix-vector products.
To deal with this term efficiently, we use the trick from~\cite{ksv-manprec-2016} and assume that 
\begin{equation}\label{eq:assump_prec}
	\Prec = \mathbf{B}_1 + \dots + \mathbf{B}_{\rho_{\mathbf{B}}},
\end{equation}
where $\mathbf{B}_i$, $i=1,\dots,\rho_{\mathbf{B}}$ are of TT-rank $1$.
This trick helps us thanks to the fact the multiplication of a TT-matrix of TT-rank 1 by a general TT-matrix does not change the rank of the latter.
The assumption~\eqref{eq:assump_prec} holds for the preconditioner we use for the calculation of vibrational spectra of molecules, while for spin chains computation no preconditioner is used.
Thus,
\begin{equation}\label{eq:trick_prec}
	\OP_k\Prec\mH \X_k = \OP_k\mathbf{B}_1\mH \X_k + \dots + \OP_k\mathbf{B}_{\rho_{\mathbf{B}}}\mH \X_k
\end{equation}
and hence the complexity of computing~$\OP_k\Prec\mH \X_k$ is $\mathcal{O}(\bsize dn^2 r^2 R^2\rho_{\mathbf{B}})$.
The truncation operation costs $\mathcal{O}(\bsize dnr^3)$, which is negligible compared to other operations.
The overall algorithm is summarized in Algorithm~\ref{alg:riemannian_lobpcg}.

\begin{algorithm}[t]
\begin{algorithmic}[1]
\Require TT-matrix $\mH$, initial guess $\bX_1 = [\mX^{(1)}_1 \dots \mX^{(\bsize)}_1]$, where $\mX^{(1)}_i$ are TT-tensors, TT-rank $r$, convergence tolerance $\varepsilon$ 
\Ensure $\bX_k = [\mX^{(1)}_k \dots \mX^{(\bsize)}_k]$
 \State Initialize $\PP_k = 0 \cdot \bX_1$
 \State Set $t_1 = 1$
 \For{$k=1, 2, \dots$ until converged} 
 	\If{$\|\OP_{\mX^{(1)}_k} \left(\mH \mX_k^{(1)} -\RQ (\mX_k^{(1)}) \mX_k^{(1)} \right) \| > \varepsilon $)} 
 		\State Set $t_k = 1$
 	\Else
 		\State Choose $t_k$ randomly or according to~\eqref{eq:schedule}
 	\EndIf
 	\State Calculate $\OP_k \bX_k$ and $\OP_k \PP_k$   \Comment{$\OP_k \equiv \OP_{\mX_k}^{(t_k)}$}
 	\State Calculate $\RQ (\mX_k^{(i)})$, $i=1,\dots,\bsize$
 	\State Calculate $\OP_k\Prec \R_k$ using \eqref{eq:assump_prec} and \eqref{eq:trick_prec}
 	\State Set $\basis_k = \OP_k \left [ \X_k, \Prec \R_k, \PP_k \right ]$ and calculate $\basis_k^\intercal \basis_k$, $\basis_k^\intercal \mH \basis_k$, $\bX_k^\intercal \bX_k$, $\basis_k^\intercal \bX_k$, $\basis_k^\intercal \mH \basis_k$, $\basis_k^\intercal \mH \X_k$, $\mathrm{diag}(\bX_k^\intercal\mH \bX_k)$ as described in Sec.~\ref{sec:inner}
 	\For {$i=1, \ldots, \bsize$}
 		\State $\Rvec_k^{(i)} = \mH \mX_k^{(i)} - \RQ (\mX_k^{(i)}) \mX_k^{(i)}$
  \EndFor
  \State $\brel, \C$ = \texttt{find\_coefficients}($\basis_k^\intercal \basis_k$, $\basis_k^\intercal \mH \basis_k$, $\bX_k^\intercal \bX_k$, $\basis_k^\intercal \bX_k$, $\basis_k^\intercal \mH \basis_k$, $\basis_k^\intercal \mH \X_k$, $\mathrm{diag}(\bX_k^\intercal\mH \bX_k)$), see Alg.~\ref{alg:riemannian_lobpcg}.
  \State $\PP_{k+1} = [\OP_k \Prec \R_k, \OP_k \PP_k]\, \C [:, \bsize:3\bsize]$
  \State Calculate $\bX_{k+1} = \X_k\, \diag(\brel) + \OP_k \X_k \C_1 + \PP_{k+1}$
  \State Calculate $\bX_{k+1} \coloneqq \Retr_r (\bX_{k+1})$
 \EndFor
 \caption{Low-Rank Riemannian Alternating Projection LOBPCG.
 }\label{alg:riemannian_lobpcg}
\end{algorithmic}
\end{algorithm}

\section{Trace minimization problem for coefficients} \label{sec:coef}

Let us rewrite the problem~\eqref{eq:coefficients-problem} omitting index $k$ for simplicity
\begin{equation}
\label{eq:coefficients-problem1}
\begin{aligned}
& \underset{\brel,\, \C}{\text{minimize}}
& & \Tr\left( 
\begin{bmatrix}
\mathrm{diag}(\brel) \\
\C	
\end{bmatrix}^\intercal
\begin{bmatrix}
\X^\intercal \mH \X & \X^\intercal \mH \basis \\
\basis^\intercal \mH \X & \basis^\intercal \mH \basis
\end{bmatrix}
\begin{bmatrix}
\mathrm{diag}(\brel) \\
\C	
\end{bmatrix}
\right)\\
& \text{subject to}
& & \begin{bmatrix}
\mathrm{diag}(\brel) \\
\C	
\end{bmatrix}^\intercal
\begin{bmatrix}
\X^\intercal  \X & \X^\intercal  \basis \\
\basis^\intercal  \X & \basis^\intercal  \basis
\end{bmatrix}
\begin{bmatrix}
\mathrm{diag}(\brel) \\
\C	
\end{bmatrix} = \I_{\bsize}
\end{aligned}
\end{equation}
Unfortunately, this problem can not be reduced to a generalized eigenvalue problem, which can be solved by means of a reliable and optimized software packages.
Therefore, to solve it we propose an iterative method.
It is derived using the Lagrange multiplier method formally applied to the minimization problem.
The Lagrange function of~\eqref{eq:coefficients-problem1} reads
\[
\begin{split}
\mathcal{L}(\brel, \C, \Mult) = 
&\frac 12 
\Tr\left( 
\begin{bmatrix}
\mathrm{diag}(\brel) \\
\C	
\end{bmatrix}^\intercal
\begin{bmatrix}
\X^\intercal \mH \X & \X^\intercal \mH \basis \\
\basis^\intercal \mH \X & \basis^\intercal \mH \basis
\end{bmatrix}
\begin{bmatrix}
\mathrm{diag}(\brel) \\
\C	
\end{bmatrix}
\right)- \\
&\frac 12 \Tr \left[\Mult\left(
\begin{bmatrix}
\mathrm{diag}(\brel) \\
\C	
\end{bmatrix}^\intercal
\begin{bmatrix}
\X^\intercal  \X & \X^\intercal  \basis \\
\basis^\intercal  \X & \basis^\intercal  \basis
\end{bmatrix}
\begin{bmatrix}
\mathrm{diag}(\brel) \\
\C	
\end{bmatrix} - \I_{\bsize} 
\right)
\right]
\end{split}
\]
Thanks to the symmetry of the constraint matrix we can consider $\Mult^\intercal = \Mult$ without the loss of generality.
Then, the gradient of the Largangian reads
\[
\begin{split}
&\nabla_{\brel}\, \mathcal{L} = \diag(\Aoo) \brel + \diag ((\Aol) \C - (\Gol) \C \Mult) \mathbf{1} 
- (\Mult \odot \Goo) \brel, \\
&\nabla_{\C}\, \mathcal{L} = (\All) \C + (\Alo) \diag (\brel) - (\Glo) \diag(\brel) \Mult - (\Gll) \C \Mult,
\end{split}
\]
where $\diag(\mathbf{A})$ denotes a diagonal matrix with the same diagonal as $\mathbf{A}$, $\mathbf{1}$ --- vector of all ones of the corresponding size and $\mathbf{A}\odot\mathbf{B}$ --- elementwise product of matrices $\mathbf{A}$,$\mathbf{B}$ of the same size.
Thus, the critical point of the Lagrangian can be found from
\begin{equation}\label{eq:lagrange_matrix}
\begin{split}
	&\diag(\Aoo) \brel + \diag ((\Aol) \C) \mathbf{1} = \diag ( (\Gol) \C \Mult) \mathbf{1}  + (\Mult \odot \Goo) \brel, \\
    & (\All) \C + (\Alo) \diag (\brel) = (\Glo) \diag(\brel) \Mult + (\Gll) \C \Mult, \\
    & (\X\, \diag(\brel) + \basis \C)^\intercal (\X\, \diag(\brel) + \basis \C) = \I_{\bsize},
\end{split}
\end{equation}
We put emphasis on the fact that the latter equation does not depend on the whole $\X^\intercal \mH \X$ which is present in~\eqref{eq:coefficients-problem1}, but only on $\mathrm{diag}(\X^\intercal \mH \X)$ instead.
This significantly reduces complexity of computations.

Equations~\eqref{eq:lagrange_matrix} can be rewritten in the following form
\begin{equation}
\label{lagrange_detailed}
\begin{split}
		\begin{bmatrix}
            (\mX^{(\alpha)})^\intercal \mH \mX^{(\alpha)} & (\mX^{(\alpha)})^\intercal \mH \basis \\
            \basis^\intercal \mH \mX^{(\alpha)} & \basis^\intercal \mH\basis
        \end{bmatrix}
        \begin{bmatrix}
            \zeta_\alpha \\ \brel_\alpha
        \end{bmatrix}
        =
        \mult_{\alpha\alpha}
        \begin{bmatrix}
            (\mX^{(\alpha)})^\intercal  \mX^{(\alpha)} & (\mX^{(\alpha)})^\intercal  \basis \\
            \basis^\intercal  \mX^{(\alpha)} & \basis^\intercal \basis
        \end{bmatrix}
        \begin{bmatrix}
            \zeta_\alpha \\ \brel_\alpha
        \end{bmatrix} 
        &
        \\
        +
        \sum_{\substack{\beta =1, \\ \beta \not = \alpha}}^{\bsize} 
        \mult_{\alpha\beta}
        \begin{bmatrix}
            (\mX^{(\alpha)})^\intercal  \mX^{(\beta)} & (\mX^{(\alpha)})^\intercal  \basis \\
            \basis^\intercal  \mX^{(\beta)} & \basis^\intercal \basis
        \end{bmatrix}
        \begin{bmatrix}
            \zeta_\beta \\ \brel_\beta 
        \end{bmatrix} 
        &
      \\
        \begin{bmatrix}
            \zeta_\alpha & \brel_\alpha^\intercal
        \end{bmatrix}
        \begin{bmatrix}
             (\mX^{(\alpha)})^\intercal \mX^{(\beta)} & (\mX^{(\alpha)})^\intercal \basis  \\ 
            \basis^\intercal \mX^{(\beta)} &  \basis^\intercal \basis
        \end{bmatrix}
        \begin{bmatrix}
            \zeta_\beta \\ \brel_\beta
        \end{bmatrix}
        =
        \delta_{\alpha\beta}
        , \quad \alpha,\beta=1,\dots,\bsize, 
        \qquad\quad
        &
\end{split}
\end{equation}
where $\delta_{\alpha\beta}$ is the Kronecker delta and $\brel = \begin{bmatrix}
	\zeta_1 & \dots & \zeta_\bsize
\end{bmatrix}^\intercal$.
For convenience we introduce the following notations:
\begin{equation}\label{eq:coef_notation}
\begin{split}
\mathbf{A}_{\alpha} =& 
 		\begin{bmatrix}
            (\mX^{(\alpha)})^\intercal \mH \mX^{(\alpha)} & (\mX^{(\alpha)})^\intercal \mH \basis \\
            \basis^\intercal \mH \mX^{(\alpha)} & \basis^\intercal \mH\basis
        \end{bmatrix},
        \\
\mathbf{G}_{\alpha\beta} =&
  		\begin{bmatrix}
            (\mX^{(\alpha)})^\intercal  \mX^{(\beta)} & (\mX^{(\alpha)})^\intercal  \basis \\
            \basis^\intercal  \mX^{(\beta)} & \basis^\intercal \basis
        \end{bmatrix},
	\\
    \mathbf{s}_{\alpha} =& 
		\begin{bmatrix}
			 \zeta_\alpha \\ 
             \brel_\alpha
		\end{bmatrix}.
\end{split}
\end{equation}
Using these notations \eqref{lagrange_detailed} reads
\begin{equation}
\label{lagrange_short}
\begin{split}
	&\mathbf{A}_{\alpha} \mathbf{s}_{\alpha} = \lambda_{\alpha\alpha} \mathbf{G}_{\alpha\alpha}\mathbf{s}_{\alpha} + \sum_{\substack{\beta =1, \\ \beta \not = \alpha}}^{\bsize} \mult_{\alpha\beta} \mathbf{G}_{\alpha\beta}\mathbf{s}_{\beta}, 
    \\
    &\mathbf{s}_{\alpha}^\intercal \mathbf{G}_{\alpha\beta}\mathbf{s}_{\beta} = \delta_{\alpha\beta}, 
    \quad 
    \alpha,\beta = 1,\dots, \bsize.
\end{split}
\end{equation}
To solve~\eqref{lagrange_short} we propose the following iterative process:
\begin{equation}
\label{eq:lagrange_iterative}
\begin{split}
	&\mathbf{A}_{\alpha} \mathbf{s}_{\alpha}^{(k+1)} = \lambda_{\alpha\alpha}^{(k+1)} \mathbf{G}_{\alpha\alpha}\,\mathbf{s}_{\alpha}^{(k+1)} + \sum_{\beta < \alpha} \mult_{\alpha\beta}^{(k+1)} \mathbf{G}_{\alpha\beta}\,\mathbf{s}_{\beta}^{(k+1)}
    + \sum_{\beta > \alpha}^{\bsize} \mult_{\alpha\beta}^{(k)} \mathbf{G}_{\alpha\beta}\,\mathbf{s}_{\beta}^{(k)}, 
    \\
    &\left(\mathbf{s}_{\alpha}^{(k+1)}\right)^\intercal \mathbf{G}_{\alpha\beta}\, \mathbf{s}_{\beta}^{(k+1)} = \delta_{\alpha\beta}, 
    \quad 
    \alpha \leq \beta,
    \\
    &\left(\mathbf{s}_{\alpha}^{(k+1)}\right)^\intercal \mathbf{G}_{\alpha\beta}\, \mathbf{s}_{\beta}^{(k)} =  0, 
    \quad 
    \alpha > \beta.
\end{split}
\end{equation}
To reduce~\eqref{eq:lagrange_iterative} to a sequence of generalized eigenvalue problems let us start with $\alpha=1$ and proceed to $\alpha = \bsize$. For $\alpha = 2,\dots, \bsize -1$ define
\[
\begin{split}
&\nullmat_{\alpha}^{(k)} = 
	\begin{bmatrix}
	\mathbf{G}_{\alpha,1}\,\mathbf{s}_{1}^{(k+1)} & \dots & \mathbf{G}_{\alpha,\alpha-1}\,\mathbf{s}_{\alpha-1}^{(k+1)} & 
		\mathbf{G}_{\alpha,\alpha+1}\,\mathbf{s}_{\alpha+1}^{(k)} & \dots & \mathbf{G}_{\alpha,\bsize}\,\mathbf{s}_{\bsize}^{(k)}	
	\end{bmatrix},\\ 
	&\nullmat_{1}^{(k)} = 
	\begin{bmatrix}
	\mathbf{G}_{1,2}\,\mathbf{s}_{2}^{(k)} & \dots & \mathbf{G}_{1,\bsize}\,\mathbf{s}_{\bsize}^{(k)}	
	\end{bmatrix},\\ 
	&\nullmat_{\bsize}^{(k)} = 
	\begin{bmatrix}
	\mathbf{G}_{\bsize,1}\,\mathbf{s}_{1}^{(k+1)} & \dots & \mathbf{G}_{\bsize,\bsize -1}\,\mathbf{s}_{\bsize-1}^{(k+1)} 
	\end{bmatrix},
	\end{split}
\]
with null spaces
\[
	\mathcal{N}^{(k)}_\alpha = \mathrm{Null}\left({\nullmat_{\alpha}^{(k)}}^\intercal\right).
\]

We also introduce $\mathbf{Q}_\alpha$ of size $(3\bsize + 1)\times \mathrm{dim}(\mathcal{N}^{(k)}_\alpha)$ whose columns form an orthonormal basis in $\mathcal{N}^{(k)}_\alpha$.
Accounting for the fact that $\mathbf{s}_{\alpha}^{(k+1)}\in \mathcal{N}^{(k)}_\alpha$, 
multiplying the first equation of~\eqref{eq:lagrange_iterative} by $\mathbf{Q}_1^\intercal$ and introducing $\mathbf{z}_{\alpha}^{(k+1)}$: $\mathbf{s}_{\alpha}^{(k+1)} = \mathbf{Q}_\alpha \mathbf{z}_{\alpha}^{(k+1)}$ we arrive at the sequence of generalized eigenvalue problems
\begin{equation}\label{eq:eigh_mult}
\begin{split}
	&\left(\mathbf{Q}_\alpha^\intercal\mathbf{A}_{\alpha} \mathbf{Q}_\alpha\right) \mathbf{z}_{\alpha}^{(k+1)} = \lambda_{\alpha\alpha}^{(k+1)} \left(\mathbf{Q}_\alpha^\intercal\mathbf{G}_{\alpha\alpha} \mathbf{Q}_\alpha\right) \mathbf{z}_{\alpha}^{(k+1)}, \\
    &\left( \mathbf{z}_{\alpha}^{(k+1)} \right)^\intercal \left(\mathbf{Q}_\alpha^\intercal\mathbf{G}_{\alpha\alpha} \mathbf{Q}_\alpha\right) \mathbf{z}_{\alpha}^{(k+1)}  = 1,
\end{split}
\end{equation}
in which we are searching for the smallest eigenvalue~$\lambda_{\alpha\alpha}^{(k+1)}$ and the corresponding  eigenvector~$\mathbf{z}_{\alpha}^{(k+1)}$.

The matrix $\mathbf{Q}_\alpha$ is calculated using SVD of $\nullmat_{\alpha}^{(k)}=U\Sigma V^\intercal$ by choosing $\bsize$ plus number of zero singular values of last columns of $U$.
The algorithm described above is summarized in Algorithm~\ref{alg:heuristic}.
Given matrices $\basis^\intercal \basis$, $\basis^\intercal \mH \basis$, $\bX^\intercal \bX$, $\basis^\intercal \bX$, $\basis^\intercal \mH \basis$, $\basis^\intercal \mH \X$, $\mathrm{diag}(\bX^\intercal\mH \bX)$ the overall complexity to find the coefficients is $\mathcal{O}(\bsize^4)$.
This is due to the fact that we calculate full eigendecomposition that costs $\mathcal{O}(\bsize^3)$ for each $\alpha=1,\dots,\bsize$.
Note that we only need one eigenvalue and one eigenvector for each $\alpha$.
This does not significantly increase the complexity of the overall algorithm (including tensor operations) for moderate $\bsize$ up to $100$. 
Arguably $\mathcal{O}(\bsize^4)$ complexity may be reduced down to $\mathcal{O}(\bsize^3)$ using low-rank update of matrix decompositions.
We, however, do not consider it here as it does not significantly influence overall performance of the method on the considered range of $\bsize$ and requires a lot of technical details to be presented.

\begin{algorithm}[H]
\begin{algorithmic}[1]
\caption{\texttt{find\_coefficients} function from Algorithm~\ref{alg:riemannian_lobpcg}}\label{alg:heuristic}
 \Require Matrices $\basis^\intercal \basis$, $\basis^\intercal \mH \basis$, $\bX^\intercal \bX$, $\basis^\intercal \bX$, $\basis^\intercal \mH \basis$, $\basis^\intercal \mH \X$, $\mathrm{diag}(\bX^\intercal\mH \bX)$, initial guesses $\mathbf{s}_{\beta}^{(1)}\in\mathbb{R}^{3\bsize + 1}$, $\beta=2,\dots,\bsize$,
 maximum number of iterations $K$
 \Ensure $\brel, \C$ that approximate the solution to~\eqref{eq:coefficients-problem1}
 \State Assemble $\mathbf{A}_{\alpha}$, $\mathbf{G}_{\alpha\beta}$, $\alpha,\beta=1,\dots,\bsize$ according to~\eqref{eq:coef_notation}
 \State $\nullmat_{1}^{(1)}=\begin{bmatrix}
	\mathbf{G}_{1,2}\,\mathbf{s}_{2}^{(1)} & \dots & \mathbf{G}_{1,\bsize}\,\mathbf{s}_{\bsize}^{(1)}	
	\end{bmatrix}$  
 \For{$k=1, 2, \ldots, K$}
 	\For{$\alpha=1, 2, \ldots, \bsize$}
 		\State Compute SVD: $\nullmat_{\alpha}^{(k)}=U\mathrm{diag(\mathbf{\sigma})} V^\intercal$ 
 		\State Set $\mathbf{Q}_\alpha \coloneqq U[:,2\bsize+1-\#\mathrm{zeros}({\bf \sigma}):3\bsize+1]$, 
 		\State Compute $\mathbf{Q}_\alpha^\intercal\mathbf{A}_{\alpha} \mathbf{Q}_\alpha$, $\mathbf{Q}_\alpha^\intercal\mathbf{G}_{\alpha\alpha} \mathbf{Q}_\alpha$
 		\State Find $\mathbf{z}_{\alpha}^{(k+1)}$ from~\eqref{eq:eigh_mult} 
 		\State Compute $\mathbf{s}_{\alpha}^{(k+1)} = \mathbf{Q}_\alpha \mathbf{z}_{\alpha}^{(k+1)}$ 
 		\If{$\alpha\not= \bsize$}
 		    \State Calculate $\nullmat_{\alpha+1}^{(k)}=\begin{bmatrix}
	\mathbf{G}_{\alpha+1,\beta} \mathbf{z}_{\beta}^{(k+1)}
	\end{bmatrix}_{\beta=1, \beta\not =\alpha+1}^\bsize$
	 \EndIf
 	\EndFor 
 	\State $\nullmat_{1}^{(k+1)}=\begin{bmatrix}
	\mathbf{G}_{1,2}\,\mathbf{s}_{2}^{(k+1)} & \dots & \mathbf{G}_{1,\bsize}\,\mathbf{s}_{\bsize}^{(k+1)}	
	\end{bmatrix}$ 	
 \EndFor
\end{algorithmic}
\end{algorithm}

\section{Numerical experiments}

In this section, we numerically assess the proposed method.  For all experiments we use NVIDIA DGX-1 station with 8 V100 GPUs (only one of which is used at a time) and Intel(R) Xeon(R) CPU E5-2698 v4 @ 2.20GHz with 80 logical cores.

\subsection{Molecule vibrational spectra} \label{sec:num_mol}

One of the applications we consider is the computation of vibrational spectra of molecules.
In particular, we consider Hamiltonians given as
\begin{equation}\label{eq:schroed}
	\mathcal{H} = -\frac{1}{2}\sum_{i=1}^d \omega_i \frac{\partial^2}{\partial q_i^2}+ V(q_1,\dots, q_d),
\end{equation}
where $V$ denotes potential energy surface (PES).
The proposed method is applicable if the Hamiltonian can be represented in the TT-format with small TT-ranks, which holds, e.g. for sum-of-product PES.
Therefore, we consider PES given in the polynomial form as is used in~\cite{avila2011using134}:
\[
\begin{split}
V(q_1,\dots,q_{d}) = \frac{1}{2} \sum_{i=1}^{d} \omega_i q_i^2 + \frac{1}{6} \sum_{i=1}^{d} \sum_{j=1}^{d} \sum_{k=1}^{d} \phi_{ijk}^{(3)} q_i q_j q_k\\
+ \frac{1}{24} \sum_{i=1}^{d} \sum_{j=1}^{d} \sum_{k=1}^{d} \sum_{l=1}^{d} \phi_{ijkl}^{(4)} q_i q_j q_k q_l.
\end{split}
\]
To discretize the problem we use the discrete variable representation (DVR) scheme on the tensor product of Hermite meshes~\cite{baye-dvr-1986}.
The complexity of assembling the discretized PES $V$ is~\cite{ro-mp-2016}
$$\mathcal{O}\left(\left(\verb nnz \left(\phi_{ijk}^{(3)}\right) + \verb nnz \left(\phi_{ijkl}^{(4)}\right) \right) d n^2 R^3 \right)$$
and is negligibly small compared with the cost of one iteration of the iterative process.
 As a preconditioner we use approximate inversion of the harmonic part of~\eqref{eq:schroed} using~\cite{khor-prec-2009} and formula~\eqref{eq:trick_prec} for efficient computations.

We calculate vibrational spectra of two molecules: acetonitrile (CH$_3$CN) and ethylene oxide (C$_2$H$_4$O) for which $d=12$ and $d=15$ correspondingly.
 The potential energy surfaces for these molecules were kindly provided by the group of Prof. Tucker Carrington.
 Table~\ref{tab:hamranks} contains TT-ranks of Hamiltonians of these two molecules with the mode order sorted in correspondence with the ascending order of $\omega_{i}$, $i=1,\dots,d$ with mode sizes chosen according to~\cite{lc-rrbpm-2014} for CH$_3$CN and  $n=15$ for all modes for C$_2$H$_4$O.
 Initial guess is chosen from the solution of the harmonic part of the Hamiltonian and can be derived analytically.
\setlength{\tabcolsep}{4pt}
\begin{table}[t!]
\caption{TT-ranks of the Hamiltonian with $\varepsilon=10^{-12}$ truncation tolerance.}
\label{tab:hamranks}
\begin{center}
\begin{tabular}{lccccccccccccccc}
\toprule
CH$_3$CN & $R_1$ & $R_2$ & $R_3$ & $R_4$ & $R_5$ & $R_6$ & $R_7$ & $R_8$ & $R_9$ & $R_{10}$ & $R_{11}$ \\
 &  5& 9& 14&21&25&26&24&18&15&8&5\\
\midrule
C$_2$H$_4$O & $R_1$ & $R_2$ & $R_3$ & $R_4$ & $R_5$ & $R_6$ & $R_7$ & $R_8$ & $R_9$ & $R_{10}$ & $R_{11}$ & $R_{12}$ & $R_{13}$ & $R_{14}$ \\
 &  5& 11& 17& 21& 23& 25& 27 & 28 & 25 & 23 & 21 & 16 & 11 & 5\\
\bottomrule
\end{tabular}
\end{center}
\end{table}
 
\begin{figure}[t]
\begin{center}
	\includegraphics[width=120mm]{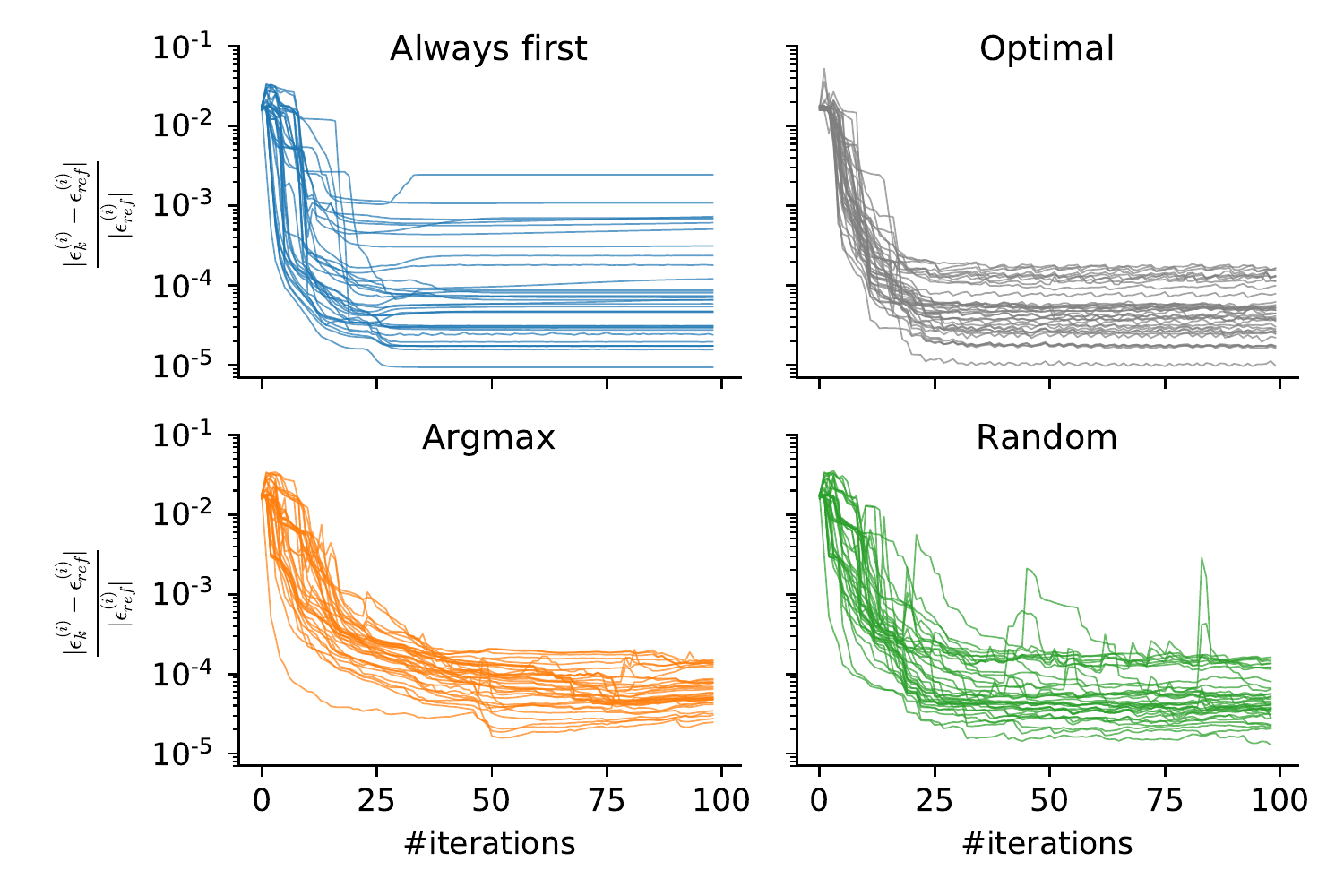}
	\caption{CH$_3$CN molecule, $r=10$. Convergence plots of $40$ eigenvalues for four different tangent space schedule scenarios: first tangent space is chosen for all the iterations, optimal choice of schedule (not practical), proposed scheduling by~\eqref{eq:schedule} and a random choice of tangent space.} \label{fig:schedule_ch3cn}
\end{center}
\end{figure}

 Figure~\ref{fig:schedule_ch3cn} represents convergence of each of the eigenvalues when running the proposed method for different choices of  tangent space schedules.
 The following scenarios are compared.
 The first one is when the tangent space is the fixed tangent space of the first eigenvector for all iterations, i.e. $t_1=\dots=t_K=1$.
In this case there is no need to find corrections using~\eqref{eq:coefficients-problem} and hence, computations are faster for large $\bsize$.
In this case we observe that although first ten eigenvalues converge within approximately $25$ iterations, most of the eigenvalues do not converge to the desired accuracy even when the number of iterations is $100$.
 The behaviour is in accordance with the fact that not necessarily can all eigenvectors be approximated using only one tangent space of the first eigenvector.
 In the optimal scenario after every iteration we choose the tangent space that allows us obtaining the smallest value of the functional in~\eqref{eq:coefficients-problem}.
 This is done by checking tangent spaces of all the eigenvectors, which is impractical and hence is provided only for comparison purposes.
 Not surprisingly, this scenario yields the fastest convergence.
 Finally, in the figure we also provide two practically interesting cases: tangent space schedule according to~\eqref{eq:schedule}, which we call here ``argmax'' and the random choice of tangent space after each iteration.
The ``argmax'' strategy aims at speeding up convergence of the eigenvalues with the slowest convergence rate, thus certifying convergence of larger eigenvalues.
By contrast, the random choice of tangent spaces provides faster convergence of the smallest eigenvalues, while the error of larger eigenvalues fluctuates.
In the following numerical experiments we choose the ``argmax'' strategy as a more reliable one and combine it with the usage of the first tangent space for the first $20$ iterations.

\begin{figure}[t]
\centering
\subfloat[]{\label{fig:time_vs_b_ch3cn}{\includegraphics[width=0.45\textwidth]{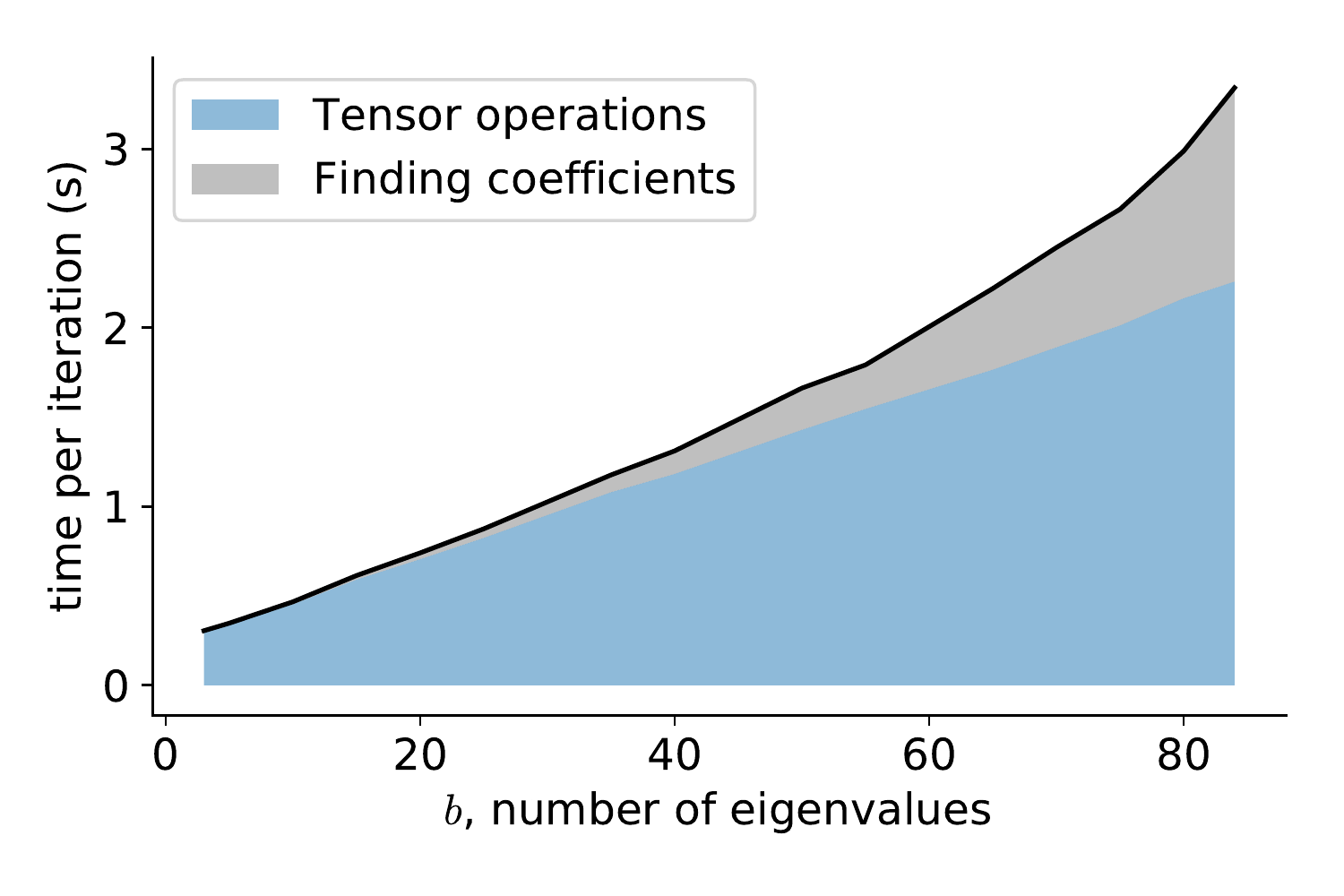}}}
\subfloat[]{\label{fig:time_vs_b_spins}{\includegraphics[width=0.45\textwidth]{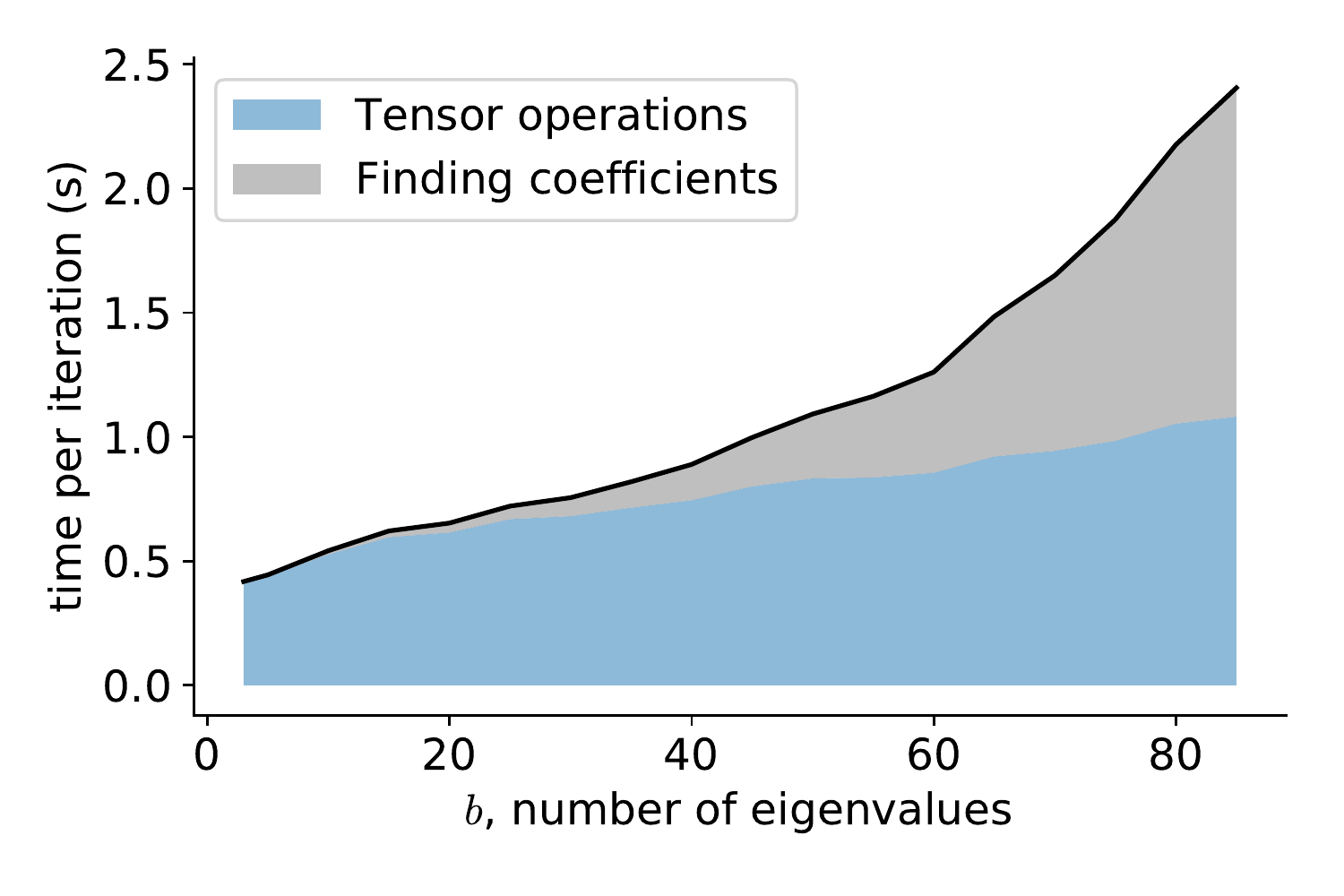}}}
\caption{GPU energy spectra calculation for acetonitrile molecule CH$_3$CN (a) and for spin lattice with $d=40$ spins (b). Time per one iteration vs. $\bsize$ -- number of computed eigenvectors, TT-rank for the both cases is $25$. Plot illustrates linear scaling of tensor operations complexity with respect to $\bsize$. For large $b$ problem of finding coefficients (Alg.~\ref{alg:heuristic}) starts dominating.} 
\label{fig:time_vs_b}
\end{figure}

In Figure~\ref{fig:time_vs_b_ch3cn} we provide computational time of one iteration with respect to the number of computed eigenvalues $\bsize$ for acetonitrile molecule CH$_3$CN and $r=25$.
This figure illustrates that the tensor part of computations described in Section~\ref{sec:impl} scales effectively linearly in the given range of $\bsize$. 
When $\bsize$ approaches $100$, the problem of finding coefficients (Alg.~\ref{alg:heuristic}) starts dominating.
We plan to improve the complexity of this part of the method in the future work.

Table~\ref{tab:errs_mols} represents accuracies of energy levels for different methods for CH$_3$CN and C$_2$H$_4$O molecules.
We measure the mean absolute error (MAE) of energy levels:
\begin{equation}\label{eq:mae}
    \mathrm{MAE} = \frac{1}{\bsize} \sum_{i=1}^{\bsize} \left|\tilde\epsilon^{(i)} - \epsilon^{(i)}_\mathrm{ref}\right|,
\end{equation}
where $\tilde\epsilon^{(i)}$, $i=1,\dots,\bsize$ are the calculated energy levels and $\epsilon^{(i)}_\mathrm{ref}$ are accurate energy levels calculated in~\cite{ro-mp-2016} for CH$_3$CN and in~\cite{intertwined2017thomas} for C$_2$H$_4$O.
For comparison purposes we also provide timings of the methods.
The timings of LRRAP LOBPCG and MP LOBPCG are measured on the same machine, whereas timings of hierarchical rank reduced block power method (H-RRBPM) are taken from~\cite{tc-hrrbpm-2015}.
We do not provide timings of H-RRBPM and HI-RRBPM for C$_2$H$_4$O since the calculations in~\cite{intertwined2017thomas} were done for different $\bsize$ ($\bsize=200$) compared with $\bsize=35$ for LRRAP LOBPCG and MP LOBPCG.
Table~\ref{tab:errs_mols} illustrates that the method is capable of producing accurate results with speedups up to 20 times on GPU compared with CPU.
We note that for larger molecules additional speedup can be obtained using several GPUs simultaneously.
In contrast to the MP LOBPCG method~\cite{ro-mp-2016} the proposed method is capable of producing comparably accurate results faster both on CPUs and GPUs.
The fact is that for the considered examples the cost of each iteration of the proposed method is considerably less than the cost of an iteration of the MP LOBPCG method, although LRRAP LOBPCG requires more iterations.
Moreover, MP LOBPCG introduces errors after such operations as calculation of linear combinations of vectors and matrix-vector products due to truncation errors, which eventually leads to lower accuracies of the result.
At the same time in the LRRAP approach corrections on each iteration  belong to a single tangent space, so no rank growth occurs.
Moreover, before the retraction all tensor calculations are done with the machine precision.
Table~\ref{tab:errs_mols} also illustrates that LRRAP LOBPCG for $r=25$ is  more accurate than the most accurate basis (basis-3 or ``b3'' for short) considered in~\cite{tc-hrrbpm-2015}.
Acceleration on GPUs allows to get additional gain in time w.r.t. H-RRBPM method.
We note, however, that the recently proposed HI-RRBPM~\cite{intertwined2017thomas} and its improved version~\cite{intertwined2018thomas} are superior to H-RRBPM.
We also note that accuracy of eigenvalues can be additionally improved using the manifold-projected simultaneous inverse iteration (MP SII) proposed for the TT-format in~\cite{ro-mp-2016}.
This strategy was used to correct eigenvalues obtained using MP LOBPCG with small~$r$.

\setlength{\tabcolsep}{4pt}
\begin{table}[t!]
\caption{Mean absolute error (MAE)~\eqref{eq:mae} and computation times of different methods.  
Reference energies $\epsilon^{(i)}_\mathrm{ref}$ are taken from~\cite{ro-mp-2016} ($r=40$) for CH$_3$CN  and from~\cite{intertwined2017thomas} (setting D) for C$_2$H$_4$O. $\bsize=84$ energy levels are calculated for CH$_3$CN and $\bsize=35$ energy levels are calculated for molecule C$_2$H$_4$O. Note that timings of H-RRBPM method are taken from~\cite{tc-hrrbpm-2015}, so computations were performed on a different machine. MAE error for H-RRBPM and HI-RRBPM is measured for first $\bsize$ eigenvalues.}
\label{tab:errs_mols}
\begin{center}
\begin{tabular}{c|cc|cc|ccc|cccc}
\toprule
  & \multicolumn{2}{c|}{\small\makecell{LRRAP \\ LOBPCG \\ (proposed)}} & \multicolumn{2}{c|}{\small\makecell{MP \\ LOBPCG \\ \cite{ro-mp-2016}}} &  \multicolumn{3}{c|}{\small\makecell{H-RRBPM\\ \cite{tc-hrrbpm-2015}}}  & \multicolumn{3}{c}{\small\makecell{HI-RRBPM \\ \cite{intertwined2017thomas}}}  \\
 \midrule
 \cellcolor{black!10}{CH$_3$CN} & {\small $r\shorteq 15$} & {\small $r\shorteq 25$} & {\small $r\shorteq 15$} & {\small $r\shorteq 25$} & b1 & b2 & b3 & \\
 \midrule
 {\small MAE,\,cm$^{-1}$} & $0.4$ & $0.05$ & $0.5$ & $0.07$ & $2.6$ & $0.9$ & $0.2$     \\
{\small GPU time} & $51$ s & $114$ s & & & & & & &  \\
{\small CPU time} & 12\,{\small min} & 26\,{\small min}  & $27$\,{\small min} & $45$\,{\small min} & $44$ s & $11$\,{\small min} & $3.2$\,{\small h} \\
 \midrule
 \cellcolor{black!10}{C$_2$H$_4$O} & $r\shorteq 25$ & $r\shorteq 35$ & $r\shorteq 25$ & $r\shorteq 35$  & \multicolumn{3}{c|}{E {\small (H-RRBPM, \cite{intertwined2017thomas})}}  & A & B & C \\
\midrule
 {\small MAE,\,cm$^{-1}$} & $1.5$ & $0.3$ & $1.7$ & $0.4$ & \multicolumn{3}{c|}{0.78} & $1.0$ & $0.6$ & $0.2$  \\
{\small GPU time} & 92 s & 129 s & & & & & & & \\
{\small CPU time} & $25$\,{\small min} & $41$\,{\small min} & $32$\,{\small min} & $64$\,{\small min} & & & & & \\
\bottomrule
\end{tabular}
\end{center}
\end{table}

 \subsection{Spin chains} \label{sec:task_spins}
As a second application, we consider Heisenberg model for one-dimensional lattices of spins.
 The goal it to compute minimal energy levels of Hamiltonian
 \begin{equation}\label{eq:ham_spin}
  \mH =	\sum_{i=1}^{d-1} \left( \bS_x^{(i)} \bS_x^{(i+1)} + \bS_y^{(i)}\bS_y^{(i+1)} + \bS_z^{(i)} \bS_z^{(i+1)} \right),
 \end{equation}
 with
 \[
 	\bS_\alpha^{(i)} = \I\otimes \dots \otimes \I \otimes \bS_\alpha \otimes \I \otimes \dots \otimes \I, \quad \alpha = x,y,z,
 \]
 where $\I$ is $2\times 2$ identity matrix and $\bS_\alpha$, $\alpha = x,y,z$ are elementary Pauli matrices
 \[
 	\bS_x = 
     \frac{1}{2}
 	\begin{bmatrix}
 	  1	& 0 \\
 	  0 & -1
 	\end{bmatrix},
 	\quad
 	\bS_y = 
     \frac{\mathrm{i}}{2}
 	\begin{bmatrix}
 	  0	& 1 \\
 	  -1 & 0
 	\end{bmatrix},
 	\quad
 	\bS_z = 
     \frac{1}{2}
 	\begin{bmatrix}
 	  0	& 1 \\
 	  1 & 0
 	\end{bmatrix}.
 \]
It can be verified numerically that TT-rank of the Hamiltonian~\eqref{eq:ham_spin} is bounded from above by $5$. 
 \begin{figure}[t]
\begin{center}
	\includegraphics[width=120mm]{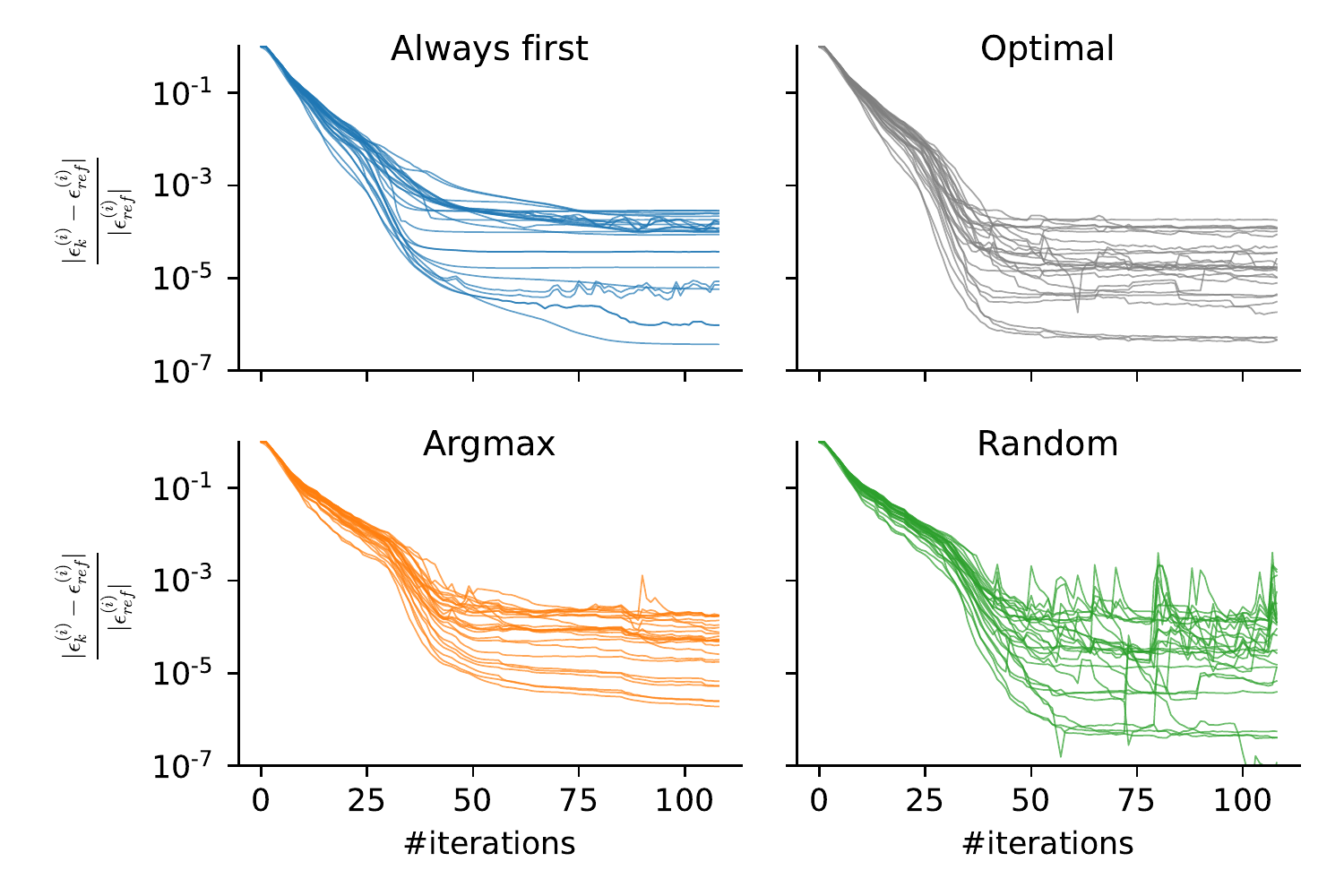}
	\caption{Heisenberg model, $\bsize=30$, $r=25$, $d=40$. Convergence plots of $30$ eigenvalues for four different tangent space schedule scenarios: first tangent space is chosen for all the iterations, optimal choice of schedule (not practical), proposed scheduling by~\eqref{eq:schedule} and a random choice of tangent space.} \label{fig:schedule_spins}
\end{center}
\end{figure}
 
 Similarly to Sec.~\ref{sec:num_mol} in Figure~\ref{fig:schedule_spins} we plot convergence of each of the eigenvalues for different choices of  tangent space schedules: when only tangent space of the first eigenvector is used, strategy with the optimal choice of schedule, ``argmax'' strategy~\eqref{eq:schedule} and a random choice of tangent spaces.
 The convergence behaviour is similar to the convergence behaviour of vibrational spectra computation (see Sec.~\ref{sec:num_mol}).
 The only difference we observe is that all eigenvectors can be well approximated in the tangent space of the first eigenvector.
 This allows to run most of the iterations in one tangent space.
 After that only a few iterations of the ``argmax'' strategy are needed to increase accuracy, which leads to a significant complexity reduction. 
 Note that by contrast to the computation of vibrational spectra no preconditioner is used for~\eqref{eq:ham_spin}.
 
 \setlength{\tabcolsep}{4pt}
\begin{table}[t!]
\caption{Mean absolute error (MAE)~\eqref{eq:mae} and computation times of different methods for open spin chains.  
Reference energies $\epsilon^{(i)}_\mathrm{ref}$ are calculated using \texttt{eigb}~\cite{dkos-eigb-2014} with $\delta=10^{-5}$, where $\delta$ denotes the truncation error for both \texttt{eigb} and ALPS.}
\label{tab:errs_spins}
\begin{center}
\begin{tabular}{c|ccc|cc|cccccc}
\toprule
  & \multicolumn{3}{c|}{\small\makecell{LRRAP \\ LOBPCG \\ (proposed)}} & \multicolumn{2}{c|}{\small\makecell{eigb \\ \cite{dkos-eigb-2014}}} &
  \multicolumn{2}{c}{\small\makecell{ALPS \\ \cite{alps2011bauer}}} \\
 \midrule
 \cellcolor{black!10}{$d=40$, $\bsize=5$} & {\small $r \shorteq 20$} & {\small $r \shorteq 35$} &  {\small $r \shorteq 45$}& {\small $\delta \shorteq 10^{-2}$} & {\small $\delta \shorteq 10^{-3}$}  & {\small $\delta \shorteq 10^{-4}$} & {\small $\delta \shorteq 10^{-5}$}  \\
 \midrule
 {\small MAE} & 1.0e-4 & 1.2e-5 & 2.2e-6 & 2.2e-4  & 2.4e-6 &   1.0e-2 & 1.2e-3    \\
{\small GPU time} & 9\,s  & 25\,s  & 44\,s &  & & &  \\
{\small CPU time} & 50\,s & 145\,s  & 251\,s & 6\,s & 26\,s  & 35\,s & 44\,s\\
\midrule
 \cellcolor{black!10}{$d=40$, $\bsize=35$} & {\small $r \shorteq 20$} & {\small $r \shorteq 35$} &  {\small $r \shorteq 45$}& {\small $\delta \shorteq 10^{-2}$} & {\small $\delta \shorteq 10^{-3}$}  & {\small $\delta \shorteq 10^{-4}$} & {\small $\delta \shorteq 10^{-5}$}  \\
 \midrule
  {\small MAE} & 7.7e-3 & 1.3e-4 & 5.1e-6 & 1.3e-4  & 3.0e-6  &   1.4e-2 & 1.6e-3   \\
{\small GPU time} & 28\,s  & 51\,s  & 85\,s &  & &  \\
{\small CPU time} & 5.6\,{\small min} & 12\,{\small min}  & 21\,{\small min} & 18\,{\small min} &31\,{\small min} &   24\,{\small min} & 41\,{\small min} \\
\bottomrule
\end{tabular}
\end{center}
\end{table}

In Figure~\ref{fig:time_vs_b_spins} time of one iteration on GPU is provided for $r=25$ and $d=40$.
By contrast to the computation of molecule vibrational energy levels (Fig.~\ref{fig:time_vs_b_ch3cn}), for $\bsize=84$ the time spent for tensor operations is less than the time to find coefficients.
This is due to the fact that TT-ranks of molecule Hamiltonians (Tab.~\ref{tab:hamranks}) are larger than those of the considered spin lattices, where TT-rank is bounded by $5$.

We compare the proposed method with two open source packages.
The first one is \texttt{eigb} method~\cite{dkos-eigb-2014} available in the~\texttt{ttpy}\footnote{\url{https://github.com/oseledets/ttpy}} library.
The second one is DMRG-based algorithm implemented in ALPS~\cite{alps2011bauer} (algorithms and libraries for physical simulations), which also allows to compute more than one of eigenstates.
We present results of the comparison in Table~\ref{tab:errs_spins}.
We observe that for small $\bsize$ ($\bsize=5$) we are not able to be faster than both~\texttt{eigb} and ALPS at comparable accuracies. 
For larger $\bsize$ ($\bsize=35$) the proposed method is faster on CPU, and GPU provides additional acceleration up to $15$ times.
We note that with our method we are capable of calculating larger $\bsize$ (up to $100$) with no problems, while~\texttt{eigb} struggles in this range.
The point is that due to the usage of the block TT format, the rank in~\texttt{eigb}  rapidly grows with $\bsize$. Therefore,  much larger rank values (more than $1000$ for $\bsize>40$~\cite{dkos-eigb-2014}) are required for~\texttt{eigb}.

\section{Related work}

The computation of energy levels of multidimensional Hamiltonians using low-rank tensor approximations has been considered in several communities.
In solid state physics computation of energy levels of one-dimensional spin lattices using tensor decompositions has been known for a long time.
For this kind of applications the \emph{matrix product state} (MPS) representation~\cite{fannes-mps-1992,klumper-mps-1993} is used to approximate eigenfunctions. 
MPS is known to be equivalent to the Tensor-Train format~\cite{osel-tt-2011} used in the current paper.
The classical algorithm to approximate ground state of a one-dimensional spin lattice using MPS representation is the \emph{density matrix renormalization group} (DMRG)~\cite{white-dmrg-1992, ostlund-dmrg-1995}, see review~\cite{schollwock-2011} for details.
Although MPS and DMRG are widely used in solid state physics, they were unknown in numerical analysis.

Calculation of several eigenvalues using two-site DMRG was first considered in papers by S.~White~\cite{white-dmrg-1992,white-dmrg-1993}.
Alternatively, one can use the \emph{numerical renormalization group} (NRG)~\cite{wilson-nrg-1975} and its improved version \emph{variational NRG}~\cite{pizorn-eigb-2012}.
In these methods index, corresponding to the number of an eigenvalue is present only in the last site (core) of the MPS representation.
The \texttt{eigb} algorithm~\cite{dkos-eigb-2014}, which is an extension of~\cite{khos-dmrg-2010}, alternately assigns eigenvalue index to different sites (cores) of the representation, which allows for rank adaptation.
By contrast, in the proposed method we consider every eigenvector separately, which implies that eigenvectors do not share sites (cores) of the decomposition.
This leads to smaller rank values of eigenvectors.
Moreover, thanks to the Riemannian optimization approach no rank growth occurs at iterations.

Calculation of energy levels using tensor decompositions was also considered for molecule vibrational spectra computations.
CP-decomposition~\cite{kolda-review-2009} (canonical decomposition) of eigenvectors of vibrational problems was considered in~\cite{lc-rrbpm-2014}, where rank reduced block power method (RRBPM) was used.
Its hierarchical version (H-RRBPM) was later proposed in~\cite{tc-hrrbpm-2015}.
The H-RRBPM was improved in~\cite{intertwined2017thomas} and later in~\cite{intertwined2018thomas} (HI-RRBPM, where ``I'' stands for ``intertwined''). 
With HI-RRBPM vibrational energy levels of molecules with more than dozens of atoms such as uracil and naphthalene~\cite{intertwined2018thomas} were accurately calculated.
We also note that tensor decompositions were used in quantum dynamics, and in particular in the multiconfiguration time dependent Hartree (MCTDH) method \cite{meyer-book-2009} as well as in its multilayer version~\cite{wt-mlmctdh-2003} (ML-MCTDH), which is similar to the Hierarchical Tucker representation~\cite{gras-hsvd-2010}.

Finally, a number of methods was developed independently in the mathematical community.
Tensor versions of power and inverse power methods (inverse iteration) were considered in~\cite{beylkin-2002,garcke-mregr-2009,khst-eigen-2012}.
For tensor versions of preconditioned eigensolvers such as preconditioned inverse iteration (PINVIT) and LOBPCG methods see~\cite{mach-innereig-2013,ro-hf-2016,ro-crossconv-2015,lebedeva-tensornd-2011,tobler-htdmrg-2011}.
The generalization of the Jacobi-Davidson method was considered in~\cite{ro-jd-2018}.
The solvers based on alternating optimization procedures such as ALS~\cite{holtz-ALS-DMRG-2012} or AMEn~\cite{ds-amen-2014} are proposed in~\cite{ds-dmrgamen-2015,kressner-evamen-2014}.

\section{Conclusion and future work}

We propose an eigensolver for high-dimensional eigenvalue problems in the TT format (MPS).
The ability of the solver to efficiently calculate up to $100$ eigenstates is assured by the usage of Riemannian optimization, which allows to avoid the rank growth naturally.
The solver is implemented in TensorFlow, thus allowing both CPU and GPU parallelization.
In the considered numerical examples the GPU version of the solver produces 15-20 times acceleration compared with the CPU version.

At each iteration of the solver, there arises a small, but nonstandard optimization problem to find coefficients of the iterative process.
The method proposed to address this problem allows to solve it with the complexity comparable to the complexity of tensor operations for up to $\bsize=100$ eigenvalues.
For $\bsize>100$ there is still room for improvement, which we plan to address in the future work.

\section*{Acknowledgements}

This work was supported by a grant of the Russian Science Foundation (Project No. 17-12-01587).

\newpage
\bibliographystyle{unsrt}
\bibliography{vbib,bibtex/tensor,bibtex/our,bibtex/dmrg,bibtex/iter,bibtex/molecular}

\appendix
\section{TensorFlow implementation overview} \label{sec:app-t3f-implementation}
In this section, we provide brief introduction into TensorFlow (which we use to simplify GPU support) and details on how we implement tensor operations.

We implemented the proposed method in Python relying on two libraries: TensorFlow and T3F.
TensorFlow~\cite{tf-2015} is a library written by Google to use for Machine Learning research (i.e., fast prototyping) and production alike. The focus is the ease of prototyping, GPU support, good parallelization abilities, and automatic differentiation\footnote{Given a computer program which defines a differentiable function, automatic differentiation generates another program which can compute the (exact) gradient of the function in the time at most 4 times slower than executing the original function. (In this paper we don't use automatic differentiation.)}.

TensorFlow provides a library of linear-algebra (and some other) functions which abstract away the hardware details. For example, running the matrix multiplication function \texttt{tensorflow.matmul(A, B)} would multiply the two matrices on a CPU using all the available threads. Running the same code on a computer with a GPU will execute the matrix multiplication on the GPU. When executing a sequence of TensorFlow operations, TensorFlow makes sure that the data is not copied to and from GPU memory in between the operations. For example, running
\begin{verbatim}
  A = numpy.random.normal(100, 100)
  matmul = tensorflow.matmul(A, A)
  print(tensorflow.reduce_sum(matmul))
\end{verbatim}
will generate a random matrix using Numpy on CPU, then copy it to the GPU memory, perform matrix multiplication on GPU, then compute the sum of all elements in the resulting matrix on GPU, and only then copy the result back to the CPU memory for printing.
Moreover, TensorFlow allows to easily process pieces of data on multiple GPUs and then combine the result together on a single master GPU.

Another important TensorFlow feature is vectorization. Almost all functions support working with \emph{batches} of objects, i.e. executing the same operations on a set of different arrays. For example, applying \texttt{tensorflow.matmul(A, B)} to two arrays of shapes $100 \times 3 \times 4$ and $100 \times 4 \times 5$ will return an array \texttt{C} of shape $100 \times 3 \times 5$ such that $C^{(i)} = A^{(i)} B^{(i)}$, $i = 1, \ldots, 100$, and all the matrices are multiplied in parallel. This is especially important on GPUs: because of massive parallel resources of GPUs, running 100 small matrix multiplications sequentially (in a for loop) is almost 100 times slower than multiplying them all in parallel in a single vectorized operation.

The second library used in this paper is Tensor Train on Tensor Flow (T3F)~\cite{novikov2018tensor}, which is written on top of TensorFlow and provides many primitives for working with Tensor Train decomposition.
To speed up the computations in the proposed method, we represent the current approximation to the eigenvectors $\X = \{\mX^{(1)}, \ldots, \mX^{(\bsize)}\}$ as a \emph{batch} of $\bsize$ TT-vectors, letting T3F and TensorFlow vectorize all the operations w.r.t. the number of TT-vectors. We also treat the basis $\basis$ on each iteration as a batch of TT-vectors. T3F library supports batch processing and for example executing \texttt{t3f.bilinear\_form(A, $\X$, $\X$)} finds the value of $(\mX^{(i)})^\intercal A \mX^{(i)}$ for each $i = 1, \ldots, \bsize$ in parallel on a GPU. When dealing with large problems, we also use the multigpu feature of TensorFlow to use all the available GPUs on a single computer.

Let us consider a detailed example of adding two TT-tensors together. Given two tensors $\tens{A},\tens{B}\in \mathbb{R}^{n \times \ldots \times n}$ represented in the TT-format
$$
\tensel{A}_{i_1 \ldots i_d} = G^\tens{A}_1(i_1) G^\tens{A}_2(i_2) \ldots G^\tens{A}_d(i_d) 
$$
$$
\tensel{B}_{i_1 \ldots i_d} = G^\tens{B}_1(i_1) G^\tens{B}_2(i_2) \ldots G^\tens{B}_d(i_d) 
$$
where for any $i_k = 1, \ldots, n$ the matrix $G^\tens{A}_k(i_k)$ is of size $r \times r$ for $k=2, \ldots, d-1$. $G^\tens{A}_1(i_1)$ is of size $1 \times r$, and $G^\tens{A}_d(i_d)$ is of size $r \times 1$.
The goal is to find TT-cores $\{G^\tens{C}_k(i_k)\}_{k=1}^d$ of tensor $\tens{C} = \tens{A} + \tens{B}$. The expression for these TT-cores is given by~\cite{osel-tt-2011}:
\begin{equation}
\begin{aligned}
\label{eq:add_tt_cores}
G^\tens{C}_k(i_k) &= 
\begin{bmatrix}
G^\tens{A}_{k}(i_{k}) & 0\\ 
0 & G^\tens{B}_{k}(i_{k})
\end{bmatrix}, ~~k = 2, \ldots, d-1\\
G^\tens{C}_1(i_1) &= 
\begin{bmatrix} G^\tens{A}_1(i_1) & G^\tens{B}_1(i_1) \end{bmatrix},
\quad 
G^\tens{C}_d(i_d) = 
\begin{bmatrix}
G^\tens{A}_d(i_d)\\ 
G^\tens{B}_d(i_d)
\end{bmatrix}
\end{aligned}
\end{equation}

Note that we also want to support batch processing, e.g. adding 100 TT-tensors $\tens{C}^{(i)} = \tens{A}^{(i)} + \tens{B}^{(i)}$, $i = 1, \ldots, 100$. The resulting program is similar to any implementation of summation of TT-tensors with two exceptions: support of batch processing and using TensorFlow for all elementary operations to allow GPU support.

T3F represents a batch of $d$-dimensional TT-tensors as a list of $d$ arrays, $k$-th array representing $k$-th TT-core of all the tensors in the batch. The shape of the $k$-th array ($k = 2, \ldots, d-1$) is $b \times r \times n \times r$, where $b$ is the batch-size, the shape of the first array ($k=1$) is $b \times 1 \times n \times r$, and the shape of the last array ($k = d$) is $b \times r \times n \times 1$.

\begin{algorithm}
\begin{algorithmic}[1]
\Require Arrays representing TT-cores of batches of tensors $\{\tens{A}^{(i)}\}_{i=1}^b$ and $\{\tens{B}^{(i)}\}_{i=1}^b$
\Ensure Array representing TT-cores of batch $\{\tens{C}^{(i)}\}_{i=1}^\bsize = \{\tens{A}^{(i)} + \tens{B}^{(i)}\}_{i=1}^\bsize$
 \State Concatenate $G^\tens{A}_1$ and $G^\tens{B}_1$ (of shape $b \times 1 \times n \times r$) along the 4-th axis to form array $G^\tens{C}_1$ of shape $b \times 1 \times n \times 2 r$
 \For{$k=2, \ldots, d-1$}
 	\State Create an array of zeros of size $b \times r \times n \times r$
 	\State Concatenate $G^\tens{A}_k$ with array of zeros along $4$-th axis into array $U$ of shape $b \times r \times n \times 2r$
 	\State Concatenate array of zeros with $G^\tens{B}_k$ along $4$-th axis into array $D$ of shape $b \times r \times n \times 2r$
 	\State Concatenate $U$ and $D$ along $2$-nd axis into $G^\tens{C}_k$ of shape $b \times 2r \times n \times 2r$.
 \EndFor
 \State Concatenate $G^\tens{A}_d$ and $G^\tens{B}_d$ (of shape $b \times r \times n \times 1$) along the 2-nd axis to form array $G^\tens{C}_d$ of shape $b \times 2r \times n \times 1$
 \caption{Implementation of adding two batches of TT-tensors in T3F (\texttt{t3f.add(A, B)}).
 }\label{alg:t3f_add}
\end{algorithmic}
\end{algorithm}

To add two TT-tensors, T3F calls TensorFlow functions to create arrays of appropriate sizes filled with zeros  and to concatenate the TT-cores of tensors $\tens{A}$ and $\tens{B}$ with each other and with zeros (see Alg.~\ref{alg:t3f_add}).
Note that a user of T3F can ignore this implementation details and just call \texttt{t3f.add(A, B)}.

\end{document}